\documentclass[11pt, a4paper]{article}
\usepackage{wasysym}
\usepackage{latexsym}
\usepackage{amsfonts}
\usepackage{mathrsfs}
\usepackage{amssymb}
\usepackage{dsfont}
\usepackage[all]{xy}

\newtheorem{Definitions1}{Definition}[section]

\newtheorem{Theorems1}{Theorem}[section]

\newtheorem{Coroll1}[Theorems1]{Corollary}
\newtheorem{Lemma1}[Theorems1]{Lemma}

\newtheorem{Quest1}{Question}[section]

\newenvironment{proof}[1][Proof]{\begin{trivlist}
\item[\hskip \labelsep {\bfseries #1}]}{\end{trivlist}}

\begin{document}
\title{Automorphisms of models of set theory and extensions of $\mathrm{NFU}$}
\author{Zachiri McKenzie}

\maketitle

\begin{abstract}
In this paper we exploit the structural properties of standard and non-standard models of set theory to produce models of set theory admitting automorphisms that are well-behaved along an initial segment of their ordinals. $\mathrm{NFU}$ is Ronald Jensen's modification of Quine's `New Foundations' set theory that allows non-sets into the domain of discourse. The axioms $\mathrm{AxCount}$, $\mathrm{AxCount}_\leq$ and $\mathrm{AxCount}_\geq$ each extend $\mathrm{NFU}$ by placing restrictions on the cardinality of a finite set of singletons relative to the cardinality of its union. Using the results about automorphisms of models of set theory we separate the consistency strengths of these three extensions of $\mathrm{NFU}$. We show that $\mathrm{NFU}+\mathrm{AxCount}$ proves the consistency of $\mathrm{NFU}+\mathrm{AxCount}_\leq$, and $\mathrm{NFU}+\mathrm{AxCount}_\leq$ proves the consistency of $\mathrm{NFU}+\mathrm{AxCount}_\geq$. 
\end{abstract}

\section[Introduction]{Introduction}

In \cite{jen69} Ronald Jensen introduces a weakening of Quine's `New Foundations' ($\mathrm{NF}$), which he calls $\mathrm{NFU}$, by allowing urelements (non-sets) into the domain of discourse. Despite the innocuous appearance of this weakening Jensen, in the same paper, shows that $\mathrm{NFU}$ is equiconsistent with a weak subsystem of $\mathrm{ZFC}$ and unlike $\mathrm{NF}$ is consistent with both the Axiom of Choice and the negation of the Axiom of Infinity. In the early nineties Randall Holmes and Robert Solovay embarked upon the project of determining the relative consistency strengths of natural extensions of $\mathrm{NFU}$ and the strengths of these extensions relative to subsystems and extensions of $\mathrm{ZFC}$. The fruits of this work can be seen in \cite{hol01}, \cite{solXX}, \cite{solXY} and \cite{ena04} which pinpoint the exact strength of a variety of natural extensions of $\mathrm{NFU}$ relative to subsystems and extensions of $\mathrm{ZFC}$.\\
\\
\indent Throughout this paper we will use $\mathrm{NFU}$ to denote the theory described by Jensen in \cite{jen69} supplemented with both the Axiom of Choice and the Axiom of Infinity. We will study three extensions of $\mathrm{NFU}$ that are obtained by adding the Axiom of Counting ($\mathrm{AxCount}$), $\mathrm{AxCount}_\leq$ and $\mathrm{AxCount}_\geq$ respectively. The first of these axioms was proposed in \cite{ros53} to facilitate induction in $\mathrm{NF}$. Both $\mathrm{AxCount}_\leq$ and $\mathrm{AxCount}_\geq$ are natural weakenings of $\mathrm{AxCount}$ that were introduced by Thomas Forster in his Ph.D. thesis \cite{for77}. The combined work of Rolland Hinnion \cite{hin75} and Jensen \cite{jen69} shows that $\mathrm{NFU}+\mathrm{AxCount}_\leq$ proves the consistency of $\mathrm{NFU}$. We improve this result by showing that $\mathrm{NFU}+\mathrm{AxCount}_\leq$ proves the consistency of $\mathrm{NFU}+\mathrm{AxCount}_\geq$. We also show that $\mathrm{NFU}+\mathrm{AxCount}$ proves the consistency of $\mathrm{NFU}+\mathrm{AxCount}_\leq$.\\          
\\
\indent The techniques developed in \cite{jen69} coupled with the observations in \cite{bof88} establish a strong link between models of $\mathrm{NFU}$ and models of subsystems of $\mathrm{ZFC}$ that admit non-trivial automorphism. In light of this connection and motivated by questions related to the strength of the theory $\mathrm{NFU}+\mathrm{AxCount}_\leq$ Randall Holmes asked the following:

\begin{Quest1} \label{Q:HolmesQuestion}
Is there model $\mathcal{M} \models \mathrm{ZFC}$ that admits an automorphism $j: \mathcal{M} \longrightarrow \mathcal{M}$ such that
\begin{itemize}
\item[(i)] $j(n) \geq n$ for all $\mathcal{M} \models n \in \omega$,
\item[(ii)] $j(\alpha) < \alpha$ for some $\mathcal{M} \models \alpha \in \omega_1$? 
\end{itemize}
\end{Quest1}

\noindent A model equipped with such an automorphism would yield a model of $\mathrm{NFU}+\mathrm{AxCount}_\leq$ in which the set of infinite cardinal numbers is countable. In Sections \ref{Sec:AutomorphismsFromStandardModels} and \ref{Sec:AutomorphismsFromNonStandardModels} of this paper we construct models of subsystems of $\mathrm{ZFC}$ equipped with automorphisms that are well-behaved along initial segment of their ordinals. In Section \ref{Sec:AutomorphismsFromStandardModels} we will show that models of set theory admitting automorphisms that move no points down along an initial segment of their ordinals can be built from standard models of set theory. This result allows us to show that every complete consistent extension of $\mathrm{ZFC}$ has a model which does not move any ordinal down. We then show in Section \ref{Sec:AutomorphismsFromNonStandardModels} that models of set theory admitting automorphism which move no natural number down but do move a recursive ordinal down can be built from non-standard $\omega$-models of set theory. This allows us to give a positive answer to Question \ref{Q:HolmesQuestion} even when $\omega_1$ is replaced by $\omega_1^{\mathrm{ck}}$.\\
\\
\indent In section \ref{Sec:NFU} we describe the set theory $\mathrm{NFU}$. We will survey work in \cite{jen69} and \cite{bof88} which shows that models $\mathrm{NFU}$ can be built from models of subsystems of $\mathrm{ZFC}$ that admit non-trivial automorphism. We also survey Holmes's adaption to $\mathrm{NFU}$ of techniques developed in \cite{hin75}. These techniques show that subsystems of $\mathrm{ZFC}$ can be interpreted in extensions of $\mathrm{NFU}$. In Sections \ref{Sec:NFUPlusAxCountGEQ} and \ref{Sec:NFUPlusAxCountLEQ} we apply the model theoretic results proved in Sections \ref{Sec:AutomorphismsFromStandardModels} and \ref{Sec:AutomorphismsFromNonStandardModels} to shed light of the strength of the theories $\mathrm{NFU}+\mathrm{AxCount}_\leq$ and $\mathrm{NFU}+\mathrm{AxCount}_\geq$. It is in these sections that we separate the consistency strengths of $\mathrm{NFU}+\mathrm{AxCount}$, $\mathrm{NFU}+\mathrm{AxCount}_\leq$ and $\mathrm{NFU}+\mathrm{AxCount}_\geq$.\\
\\
\indent The study of models admitting automorphisms has also recently yielded model theoretic characterisations of set theories and other foundational theories. In \cite{ena04} Ali Enayat provides an elegant characterisation of a large cardinal extension of $\mathrm{ZFC}$ in terms of the existence of a model of a weak subsystem of $\mathrm{ZFC}$ that admits a well-behaved automorphism. Enayat's work \cite{ena06} proves similar characterisations for a variety of subsystems of second order arithmetic. We have intentionally organised the background and results relating to automorphisms of models of subsystems of $\mathrm{ZFC}$ into Sections \ref{Sec:Background}, \ref{Sec:AutomorphismsFromStandardModels} and \ref{Sec:AutomorphismsFromNonStandardModels} so readers who are only interested in these results can skip Sections \ref{Sec:NFU}, \ref{Sec:NFUPlusAxCountGEQ} and \ref{Sec:NFUPlusAxCountLEQ}.

\section[Background]{Background} \label{Sec:Background}

Throughout this article we will use $\mathcal{L}$ to denote the language of set theory. If $A, B, C, \ldots$ are new relation, function or constant symbols then we will use $\mathcal{L}_{A, B, C, \ldots}$ to denote the language obtained by adding $A, B, C, \ldots$ to $\mathcal{L}$.\\
\\
\indent Given an extension of the language of set theory $\mathcal{L}^\prime$, we will often have cause to consider the L\'{e}vy hierarchy of $\mathcal{L}^\prime$ formulae which we will denote $\Delta_0(\mathcal{L}^\prime)$, $\Sigma_1(\mathcal{L}^\prime), \Pi_1(\mathcal{L}^\prime), \ldots$. If $T$ is an $\mathcal{L}^\prime$-theory then we say that an $\mathcal{L}^\prime$-formula $\phi$ is $\Delta_n^T$ if and only if $\phi$ is provably equivalent in $T$ to both a $\Sigma_n(\mathcal{L}^\prime)$ and a $\Pi_n(\mathcal{L}^\prime)$ formula. In addition to the L\'{e}vy classes we will also have cause to consider a hierarchy of subclasses of $\mathcal{L}$ introduced by Moto-o Takahashi in \cite{tak72}. Following \cite{fk91} we define the class of $\Delta_0^\mathcal{P}$ formulae to be the class of well-formed formulae built up inductively from atomic formulae in the form $x \in y$ and bounded quantification of the form $\forall x \in y$, $\exists x \in y$, $\forall x \subseteq y$ and $\exists x \subseteq y$ using the connectives $\land$, $\lor$, $\neg$ and $\Rightarrow$. The class of all $\Sigma_{n+1}^\mathcal{P}$ formulae is defined inductively to be the class of all formulae in the form $\exists \vec{x} \psi$ where $\psi$ is $\Pi_n^\mathcal{P}$. And, the class of all $\Pi_{n+1}^\mathcal{P}$ formulae is similarly defined to be the class of all formulae in the form $\forall \vec{x} \psi$ where $\psi$ is $\Sigma_n^\mathcal{P}$.\\
\\
\indent Let $\mathcal{L}^\prime$ be a language. We will denote $\mathcal{L}^\prime$-structures, which we will also call models, using the calligraphic roman letters ($\mathcal{M}$, $\mathcal{N}$, etc.). If $\mathcal{M}$ is an $\mathcal{L}^\prime$-structure and $A$ is a relation, function or constant symbol in $\mathcal{L}^\prime$ then we will use $M^\mathcal{M}$ to denote the underlying set of $\mathcal{M}$ and $A^\mathcal{M}$ to denote the interpretation of $A$ in $\mathcal{M}$; we will write $\mathcal{M}= \langle M^\mathcal{M}, \ldots, A^\mathcal{M}, \ldots \rangle$. If $\mathcal{M}$ is an $\mathcal{L}^\prime$-structure and $\phi$ is an $\mathcal{L}^\prime$ formula or term then we write $\phi^\mathcal{M}$ for the relativisation $\phi$ to $\mathcal{M}$. If $\mathcal{M}$ is an $\mathcal{L}^\prime$-structure and $a \in M^\mathcal{M}$ then we define
$$a^*= \{ x \in M^\mathcal{M} \mid \mathcal{M} \models x \in a \}.$$
Given $\mathcal{L}^{\prime\prime} \subseteq \mathcal{L}^\prime$ and an $\mathcal{L}^\prime$-structure $\mathcal{M}$, we will use $\mathcal{M} \mid_{\mathcal{L}^{\prime\prime}}$ to denote the $\mathcal{L}^{\prime\prime}$-reduct of $\mathcal{M}$ and $\mathbf{Th}_{\mathcal{L}^{\prime\prime}}(\mathcal{M})$ to denote the $\mathcal{L}^{\prime\prime}$-theory of $\mathcal{M}$.\\
\\
\indent If $\mathcal{L}_{A, B, \ldots}$ is an extension of the language of set theory then we will use $\mathcal{L}_{A, B, \ldots}^*$ to denote the Skolemisation of $\mathcal{L}_{A, B, \ldots}$. That is, for every $\mathcal{L}_{A, B, \ldots}^*$-formula $\phi(x_1, \ldots, x_{n+1})$ there is an $n$-placed function symbol $\mathbf{t}_\phi$ in $\mathcal{L}_{A, B, \ldots}^*$. We will use $\mathbf{Sk}(\mathcal{L}_{A, B, \ldots})$ to denote the $\mathcal{L}_{A, B, \ldots}^*$-theory that includes the axioms
$$\forall \vec{y}(\exists x \phi(x, \vec{y}) \Rightarrow \phi(\mathbf{t}_\phi(\vec{y}), \vec{y})) \textrm{ for every } \mathcal{L}_{A, B, \ldots}^* \textrm{-formula } \phi(x, \vec{y}).$$
It is well known that any $\mathcal{L}_{A, B, \ldots}$-structure $\mathcal{M}$ can be expanded to an $\mathcal{L}_{A, B, \ldots}^*$-structure $\mathcal{M}^\prime$ such that $\mathcal{M}^\prime \mid_{\mathcal{L}_{A, B, \ldots}} = \mathcal{M}$ and $\mathcal{M}^\prime \models \mathbf{Sk}(\mathcal{L}_{A, B, \ldots})$. If $\mathcal{M}$ is an $\mathcal{L}_{A, B, \ldots}^*$-structure and $X \subseteq M^\mathcal{M}$ then we can build a structure $\mathcal{H}_{\mathcal{L}_{A, B, \ldots}}^{\mathcal{M}}(X)$, by inductively closing under the function symbols in $\mathcal{L}_{A, B, \ldots}^*$, with the property that $|\mathcal{H}_{\mathcal{L}_{A, B, \ldots}}^{\mathcal{M}}(X)|= |X| \cdot \aleph_0$. We call $\mathcal{H}_{\mathcal{L}_{A, B, \ldots}}^{\mathcal{M}}(X)$ the Skolem hull (of $\mathcal{M}$) generated by $X$. If $\mathcal{M} \models \mathbf{Sk}(\mathcal{L}_{A, B, \ldots})$ then $\mathcal{H}_{\mathcal{L}_{A, B, \ldots}}^{\mathcal{M}}(X) \prec \mathcal{M}$.

\begin{Definitions1}
Let $\mathcal{L}^\prime$ be a language and let $\mathcal{M}$ be an $\mathcal{L}^\prime$-structure. Let $\langle I, < \rangle$ be a linear order. We say that $\{ c_i \mid i \in I \} \subseteq M^\mathcal{M}$ is a class of order indiscernibles if and only if for all $\mathcal{L}^\prime$-formulae, $\phi(x_1, \ldots, x_n)$, and for all $i_1 < \ldots < i_n$ and $j_1 < \ldots < j_n$ in $I$,
$$\mathcal{M} \models \phi(c_{i_1}, \ldots, c_{i_n}) \iff \phi(c_{j_1}, \ldots, c_{j_n}).$$  
\end{Definitions1}

In \cite{ehr56} Ehrenfeucht and Mostowski use Ramsey's Theorem \cite{ram30} to show that every theory with an infinite model has a model that with an infinite class of order indiscernibles indexed by $\mathbb{Z}$. We will make use of the infinite version of Ramsey's Theorem in the constructions appearing in the next two sections.

\begin{Theorems1}
(Ramsey \cite{ram30}) Let $n, k \in \omega$. For all $f: [\omega]^n \longrightarrow k$, there exists an $H \subseteq \omega$ and $i \in k$ such that $|H|= \aleph_0$ and for all $x \in [H]^n$, $f(x)= i$. \Square
\end{Theorems1}

Equipped with a model that is endowed with a $\mathbb{Z}$-index class of order indiscernibles Ehrenfeucht and Mostowski observe that the Skolem hull generated by the class of indiscernibles admits a non-trivial automorphism.

\begin{Lemma1}
(Ehrenfeucht-Mostowski \cite{ehr56}) Let $\mathcal{L}^\prime$ be a language and let $\mathcal{M}$ be an $\mathcal{L}^\prime$-structure. Let $\langle I, < \rangle$ be a linear order and let $\sigma: I \longrightarrow I$ be an order automorphism. If $X= \{ c_i \mid i \in I \} \subseteq M^\mathcal{M}$ is a class of order indiscernibles then there is an automorphism $j:\mathcal{H}_{\mathcal{L}^\prime}^\mathcal{M}(X) \longrightarrow \mathcal{H}_{\mathcal{L}^\prime}^\mathcal{M}(X)$ such that for all $i \in I$, $j(c_i)= c_{\sigma(i)}$. \Square
\end{Lemma1}  

\noindent This yields the following result:

\begin{Theorems1} \label{Th:EhrenfeuchtMostowskiTheorem}
(Ehrenfeucht-Mostowski \cite{ehr56}) Let $\mathcal{L}^\prime$ be a language. If $T$ is an $\mathcal{L}^\prime$-theory with an infinite model that there exists an $\mathcal{L}^\prime$-structure $\mathcal{M} \models T$ that admits a non-trivial automorphism $j: \mathcal{M} \longrightarrow \mathcal{M}$. \Square 
\end{Theorems1}

When speaking about subsystems of $\mathrm{ZFC}$ we will refer to subschemes and extension of the separation and collection axiom schemes. If $\Gamma$ is a set of formulae then we use $\Gamma$-separation (-collection) to denote the scheme axioms that asserts separation (collection) for all formulae in $\Gamma$. We use $\mathrm{TCo}$ to abbreviate the axiom of transitive containment which says that every set is contained in a transitive set. In $\mathrm{ZFC}$ the axiom of foundation implies that every $\mathcal{L}$-definable class has an $\in$-minimal element; we call this consequence class foundation. Without transitive containment and full separation this implication breaks down. In light of this, if $\Gamma$ is a set of formulae then we use $\Gamma$-foundation to denote the scheme of axioms that asserts foundation for every class definable by a formula in $\Gamma$. Throughout this paper we will make reference to the following important subsystem of $\mathrm{ZFC}$:
\begin{itemize}
\item $\mathrm{ZFC}^-$ is $\mathrm{ZFC}$, axiomatised with collection rather than replacement, minus the powerset axiom.
\item Mac Lane set theory ($\mathrm{Mac}$) is $\mathcal{L}$-theory axiomatised by the axioms of extensionality, emptyset, pairing, union, powerset, infinity, transitive containment, $\Delta_0(\mathcal{L})$-separation and foundation. 
\item Kripke-Platek set theory ($\mathrm{KP}$) is the $\mathcal{L}$-theory obtained by deleting infinity and powerset from $\mathrm{Mac}$ and adding $\Delta_0(\mathcal{L})$-collection and $\Pi_1(\mathcal{L})$-foundation.
\item Zermelo set theory ($\mathrm{Z}$) is the $\mathcal{L}$-theory obtained by deleting transitive containment from $\mathrm{Mac}$ and adding full separation.
\item $\mathrm{KP}^{\mathcal{P}}$ is obtained by adding infinity, powerset, $\Delta_0^{\mathcal{P}}$-collection and $\Pi_1^{\mathcal{P}}$-foundation to $\mathrm{KP}$.
\end{itemize}
The theories $\mathrm{Mac}$, $\mathrm{KP}$, $\mathrm{Z}$ and $\mathrm{KP}^\mathcal{P}$ are studied extensively in \cite{mat01} which compares the strength of these systems to a variety of other subsystems of $\mathrm{ZFC}$.\\
\\ 
\indent Just as $\Delta_0(\mathcal{L})$-formulae are absolute between transitive classes, $\Delta_0^\mathcal{P}$-formulae are absolute between transitive classes that have the same notion of powerset. The following lemma is a straightforward modification of a result proved in \cite{fk91}:

\begin{Lemma1}
Let $\phi(x_1, \ldots, x_n)$ be an $\Delta_0^\mathcal{P}$-formula. If $M$ and $N$ are transitive classes such that
\begin{itemize}
\item[(i)] $M \subseteq N$, 
\item[(ii)] $\langle M, \in \rangle \models \mathrm{Mac}$, 
\item[(iii)] for all $x \in M$, 
$$\langle N, \in \rangle \models \mathcal{P}(x) \textrm{ exists},$$
\item[(iv)] for all $x \in M$, $\mathcal{P}^M(x)=\mathcal{P}^N(x)$, 
\end{itemize}
then for all $a_1, \ldots, a_n \in M$,
$$\langle M, \in \rangle \models \phi(a_1, \ldots, a_n) \textrm{ if and only if } \langle N, \in \rangle \models \phi(a_1, \ldots, a_n).$$
\Square
\end{Lemma1}
 
Let $\mathcal{L}^\prime$ be a countable extension of the language of set theory. We will write $\mathcal{L}_{\omega_1 \omega}^\prime$ for the infinitary language extending $\mathcal{L}^\prime$ that includes formulae built using conjunctions and disjunctions of countable length. We can fix a coding of the language $\mathcal{L}_{\omega_1 \omega}^\prime$ in $H_{\aleph_1}$ that is definable by a $\Delta_1^{\mathrm{KP}}$ formula over $H_{\aleph_1}$. If $\phi$ is an $\mathcal{L}_{\omega_1 \omega}^\prime$-formula then we will write $\ulcorner \phi \urcorner$ for the element of $H_{\aleph_1}$ that codes $\phi$. We say that a set $A$ is admissible if and only if it is a transitive model of $\mathrm{KP}$. If $A$ is a countable admissible set then we define
$$(\mathcal{L}_{\omega_1 \omega}^\prime)_A= \{ \phi \in \mathcal{L}_{\omega_1 \omega}^\prime \mid \ulcorner \phi \urcorner \in A \}.$$
The deduction rules of first-order logic can be extended to fragments of the infinitary logic $\mathcal{L}_{\omega_1, \omega}^\prime$ --- we will write $\Gamma \vdash_{(\mathcal{L}_{\omega_1, \omega}^\prime)_A} \sigma$ if the $\sigma$ is provable from $\Gamma$ in $(\mathcal{L}_{\omega_1, \omega}^\prime)_A$. As usual we will use $\omega_1^{\mathrm{ck}}$ to denote the Church-Kleene ordinal. If $\alpha$ is an ordinal the we will use $L_\alpha$ to denote the $\alpha^{\mathrm{th}}$ level of G\"{o}del's constructible hierarchy. It is well known that $L_{\omega_1^{\mathrm{ck}}}$ is an admissible set. Moreover, $\omega_1^{\mathrm{ck}}$ is the least $\alpha > \omega$ such that $L_\alpha$ is admissible. In \cite{bar67} Jon Barwise proves analogues of the compactness and completeness fragments of $\mathcal{L}_{\omega_1, \omega}^\prime$ coded in admissible sets.

\begin{Theorems1} \label{Th:BarwiseCompleteness}
(Barwise Completeness Theorem) Let $\mathcal{L}^\prime$ be a countable extension of the language of set theory. If $A$ is a countable admissible set then
$$\{ \ulcorner \phi \urcorner \mid (\phi \in (\mathcal{L}^\prime_{\omega_1\omega})_A) \land (\vdash_{(\mathcal{L}^\prime_{\omega_1\omega})_A} \phi) \}$$
is definable by a $\Sigma_1(\mathcal{L})$ formula over $A$. 
\Square
\end{Theorems1}

\begin{Theorems1} \label{Th:BarwiseCompactness}
(Barwise Compactness Theorem) Let $\mathcal{L}^\prime$ be a countable extension of the language of set theory. Let $A$ be a countable admissible set and let $T \subseteq (\mathcal{L}^\prime_{\omega_1\omega})_A$ be definable by a $\Sigma_1(\mathcal{L})$ formula over $A$. If for every $T^\prime \subseteq T$ with $T^\prime \in A$, $T^\prime$ has a model, then $T$ has a model.  
\Square
\end{Theorems1}

\noindent We will use the Barwise Compactness Theorem in sections \ref{Sec:AutomorphismsFromNonStandardModels} and \ref{Sec:NFUPlusAxCountLEQ} to build non-standard $\omega$-models of subsystems of $\mathrm{ZFC}$.

\section[Automorphisms from standard models of set theory]{Automorphisms from standard models of set theory} \label{Sec:AutomorphismsFromStandardModels}

In this section we will exploit the structural properties of models of set theory that are standard along an initial segment of the ordinals to produce models of set theory admitting automorphisms that are well-behaved along these initial segments. We begin by showing that the existence of an $\omega$-model of Mac Lane set theory implies the existence of an elementarily equivalent model admitting an automorphism that does not move any natural number down. Our models admitting automorphism will be built using a modification of the Ehrenfeucht-Mostowski method \cite{ehr56} using ideas from \cite{kor94}. The aim will be to build a fully Skolemised model with a $\mathbb{Z}$-indexed class of order indiscernible natural numbers. Special care will be taken in the construction of this model to ensure the following behaviour:
\begin{itemize}
\item[(i)] the class of order indiscernibles will be cofinal in the natural numbers of the Skolem hull that includes all of the order indiscernibles, 
\item[(ii)] if a natural number sitting below all of the order indiscernibles corresponds to Skolem term $\Phi$ mentioning finitely many order indiscernibles then the same natural number also corresponds to the Skolem term $\Phi^\prime$ obtained by replacing the order indiscernibles in $\Phi$ with an order equivalent set of order indiscernibles.  
\end{itemize}
The Skolem hull of the resulting model which includes the $\mathbb{Z}$-indexed class of order indiscernibles admits an automorphism which moves every natural number sitting above one of the order indiscernibles up, and fixes every natural number that sits below all of the order indiscernibles. Therefore the cofinality of the class of order indiscernibles guarantees that no natural number is moved down.\\ 
\\
\indent We begin our exposition of the results of this section by introducing an extension of the language of set theory.

\begin{Definitions1}
We define $\mathcal{L}^{\mathbf{ind}} \supseteq \mathcal{L}^*$ by adding
\begin{itemize}
\item[(i)] a unary relation $\mathcal{C}$,
\item[(ii)] constant symbols $c_i$ for each $i \in \mathbb{Z}$.
\end{itemize}
\end{Definitions1}

\noindent This language allows us to define the theory that will yield models of $\mathrm{Mac}$ admitting the desired automorphism.

\begin{Definitions1} \label{Df:TheoryW1}
We define the $\mathcal{L}^{\mathbf{ind}}$-theory $\mathbf{W}_1 \supseteq \mathrm{Mac} \cup \mathbf{Sk}(\mathcal{L})$ by adding the axioms
\begin{itemize}
\item[(i)] $\forall x (\mathcal{C}(x) \Rightarrow x \in \omega)$
\item[(ii)] $(\forall x \in \omega)(\exists y \in \omega)(y > x \land \mathcal{C}(y))$,
\item[(iii)] $\mathcal{C}(c_i)$ for all $i \in \mathbb{Z}$,
\item[(iv)] $c_i < c_j$ for all $i, j \in \mathbb{Z}$ with $i < j$,
\item[(v)] for all $\mathcal{L}^*$-formulae $\phi(x_1, \ldots, x_n)$,
$$\forall x_1 \ldots \forall x_n \forall y_1 \ldots \forall y_n \left( 
\begin{array}{c}
x_1 < \ldots < x_n \land y_1 < \ldots < y_n \land \bigwedge_{1 \leq i \leq n} (\mathcal{C}(x_i) \land \mathcal{C}(y_i))\\
\Rightarrow (\phi(x_1, \ldots, x_n) \iff \phi(y_1, \ldots, y_n))
\end{array} \right),$$
\item[(vi)] for all $\mathcal{L}^*$ Skolem functions $\mathbf{t}(x_1, \ldots, x_n)$ and for all $i_1 < \ldots < i_n$ and $j_1 < \ldots < j_n$ in $\mathbb{Z}$,
$$\forall x (\mathcal{C}(x) \Rightarrow \mathbf{t}(c_{i_1}, \ldots, c_{i_n}) < x) \Rightarrow (\mathbf{t}(c_{i_1}, \ldots, c_{i_n})= \mathbf{t}(c_{j_1}, \ldots, c_{j_n})).$$   
\end{itemize}
\end{Definitions1}

It should be noted that (i) and (ii) of Definition \ref{Df:TheoryW1} are single axioms and (iii), (iv), (v) and (vi) of Definition \ref{Df:TheoryW1} are axiom schemes. Axioms (i), (ii) and scheme (v) in Definition \ref{Df:TheoryW1} ensure that $\mathcal{C}$ is a class of indiscernibles that is cofinal in the natural numbers. We will see that the cofinality of $\mathcal{C}$ will ensure that the $c_i$s (a subclass of $\mathcal{C}$ by Definition \ref{Df:TheoryW1}(iii)) will be a cofinal subclass of the natural numbers in the Skolem hull that includes the $c_i$s. This idea of using a cofinal super class of indiscernibles to ensure that the $\mathbb{Z}$-indexed class of indiscernibles is cofinal in the Skolem hull first appears in \cite{kor94}. Axiom scheme (vi) in Definition \ref{Df:TheoryW1} ensures that the Skolem terms sent below $\mathcal{C}$ are fixed when the sequence of $c_i$s appearing in the Skolem term are replaced with an order equivalent sequence. This combined with the cofinality of the $c_i$s will ensure that the automorphism generated by the $c_i$s moves all natural numbers in one direction.\\
\\
\indent We are now able to show that the existence of an $\omega$-model of $\mathrm{Mac}$ implies the existence of an elementarily equivalent model of $\mathbf{W}_1$.

\begin{Lemma1} \label{Th:ConsistencyOfW1}
If $\mathcal{M}= \langle M^\mathcal{M}, \in^\mathcal{M} \rangle$ is an $\omega$-model of $\mathrm{Mac}$ then $\mathbf{W}_1 \cup \mathbf{Th}_{\mathcal{L}}(\mathcal{M})$ is consistent.
\end{Lemma1}

\begin{proof}
Let $\mathcal{M}= \langle M^\mathcal{M}, \in^\mathcal{M} \rangle$ be an $\omega$-model of $\mathrm{Mac}$. Let $g: \omega \longrightarrow (\omega^\mathcal{M})^*$ be an isomorphism guaranteed by the fact that $\mathcal{M}$ is an $\omega$-model. Using the Axiom of Choice in the ambient theory we can find interpretations for the Skolem functions and expand $\mathcal{M}$ to an $\mathcal{L}^*$-structure $\mathcal{M}^\prime \models \mathbf{Sk}(\mathcal{L})$ with $\mathcal{M}^\prime \mid_{\mathcal{L}}= \mathcal{M}$. We will use the Compactness Theorem to show that $\mathbf{W}_1 \cup \mathbf{Th}_{\mathcal{L}}(\mathcal{M})$ is consistent. Let $\Delta \subseteq \mathbf{W}_1$ be finite. Our first task is to choose a suitable interpretation for the predicate $\mathcal{C}$. Suppose that $\Delta$ mentions instances of the scheme (v) of Definition \ref{Df:TheoryW1} for the $\mathcal{L}^*$-formulae $\phi_0, \ldots, \phi_{m-1}$. Without loss of generality we can assume that for all $0 \leq i < m$, $\phi_i$ has arity $k$. Define $H_1: [\omega]^k \longrightarrow \mathcal{P}(m)$ by
$$H_1(\{x_1, \ldots, x_k \})= \{ i \in m \mid \mathcal{M}^\prime \models \phi_i(g(x_1), \ldots, g(x_k))\} \textrm{ where } x_1 < \ldots < x_k.$$
Now, by Ramsey's Theorem there is an infinite $C^\prime \subseteq \omega$ such that $H_1``[C^\prime]^k= \{A\}$. Let $C= g``C^\prime$. Now, $C \subseteq (\omega^\mathcal{M})^*$ is a set of elements that are order indiscernible with respect to the formulae $\phi_0, \ldots, \phi_{m-1}$. And $C$ is cofinal in $(\omega^\mathcal{M})^*$, since $\langle (\omega^\mathcal{M})^*, \in^\mathcal{M} \rangle$ is isomorphic to $\langle \omega, \in \rangle$. Therefore by interpreting $\mathcal{C}$ using $C$ we can expand $\mathcal{M}^\prime$ to a structure $\mathcal{M}^{\prime\prime}$ that satisfies (i) and (ii) of Definition \ref{Df:TheoryW1} and all instances of Definition \ref{Df:TheoryW1}(v) that are mentioned in $\Delta$.\\ 
Our task now turns to finding an interpretation for the $c_i$s that are mentioned in $\Delta$. Without loss of generality we can assume that $\Delta$ mentions exactly the constants $c_{-n}, \ldots, c_{n}$. Let $z$ be the least element of $C^\prime$. Now, suppose that $\Delta$ mentions instances of the scheme Definition \ref{Df:TheoryW1}(vi) for the Skolem terms $\mathbf{t}_0, \ldots, \mathbf{t}_{l-1}$. Without loss of generality we can assume that for all $0 \leq i < l$, $\mathbf{t}_i$ is function symbol with arity $v$. Define $H_2: [C^\prime]^v \longrightarrow (z+1)^l$ by
$$H_2(\{ x_1, \ldots, x_v \})= \langle y_0, \ldots, y_{l-1} \rangle \textrm{ where } x_1 < \ldots < x_v \textrm{ and}$$
$$y_i= \left\{
\begin{array}{ll}
\mathbf{t}_i(g(x_1), \ldots, g(x_v)) & \textrm{if } \mathcal{M}^{\prime\prime} \models \left(\begin{array}{c}
(\mathbf{t}_i(g(x_1), \ldots, g(x_v)) \in \omega) \land \\
(\mathbf{t}_i(g(x_1), \ldots, g(x_v)) < g(z))
\end{array}\right)\\
z & \textrm{otherwise}
\end{array} \right).$$
Using Ramsey's Theorem we can find an infinite $D^\prime \subseteq C^\prime$ such that $H_2``[D^\prime]^v = \{\eta\}$. Let $D= g``D^\prime$.
By interpreting $c_{(i-1)-n}$ by the $i^{\mathrm{th}}$ element of $D$ for $1 \leq i \leq 2n+1$ we expand $\mathcal{M}^{\prime\prime}$ to a structure $\mathcal{M}^{\prime\prime\prime}$ such that
$$\mathcal{M}^{\prime\prime\prime} \models \Delta \cup \mathbf{Th}_{\mathcal{L}}(\mathcal{M}) \cup \mathbf{Sk}(\mathcal{L}).$$
Therefore by compactness $\mathbf{W}_1 \cup \mathbf{Th}_{\mathcal{L}}(\mathcal{M})$ is consistent. 
\Square
\end{proof}

We now show that a model of $\mathbf{W}_1$ yields a model of $\mathrm{Mac}$ admitting an automorphism that does not move any natural number down. Lemma \ref{Th:IndiscerniblesCofinalInSkolemHull} is due to Friederike K\"{o}rner in \cite{kor94}.

\begin{Lemma1} \label{Th:IndiscerniblesCofinalInSkolemHull}
(K\"{o}rner \cite{kor94}) Let $\mathcal{M} \models \mathbf{W}_1$. For all $\mathcal{L}^*$ Skolem functions $\mathbf{t}(x_1, \ldots, x_n)$ and for all $i_1 < \ldots < i_n$ in $\mathbb{Z}$, if $\mathcal{M} \models \mathbf{t}(c_{i_1}, \ldots, c_{i_n}) \in \omega$ then there exists a $k \in \mathbb{Z}$ such that
$$\mathcal{M} \models \mathbf{t}(c_{i_1}, \ldots, c_{i_n}) < c_k.$$  
\end{Lemma1}

\begin{proof}
Let $\mathcal{M} \models \mathbf{W}_1$. Let $\mathbf{t}(x_1, \ldots, x_n)$ be an $\mathcal{L}^*$ Skolem function. Let $i_1 < \ldots < i_n$ in $\mathbb{Z}$ and assume that $\mathcal{M} \models \mathbf{t}(c_{i_1}, \ldots, c_{i_n}) \in \omega$. By Definition \ref{Df:TheoryW1}(ii) there is an $a \in (\omega^\mathcal{M})^*$ such that 
$$\mathcal{M} \models \left( \bigwedge_{1 \leq j \leq n} (c_{i_j} < a) \right) \land (\mathbf{t}(c_{i_1}, \ldots, c_{i_n}) < a) \land \mathcal{C}(a).$$
Let $k > i_n$. By indiscernibility (Definition \ref{Df:TheoryW1}(vii)):
$$\mathcal{M} \models \mathbf{t}(c_{i_1}, \ldots, c_{i_n}) < c_k.$$ 
\Square
\end{proof}

\begin{Lemma1} \label{Th:SkolemFunctionsFixedBelow}
Let $\mathcal{M} \models \mathbf{W}_1$. For all $\mathcal{L}^*$ Skolem functions $\mathbf{t}(x_1, \ldots, x_n)$ and for all $i_1 < \ldots < i_n$ and $j_1 < \ldots < j_n$ in $\mathbb{Z}$, if for all $k \in \mathbb{Z}$, $\mathcal{M} \models \mathbf{t}(c_{i_1}, \ldots, c_{i_n}) < c_k$ then 
$$\mathcal{M} \models \mathbf{t}(c_{i_1}, \ldots, c_{i_n})= \mathbf{t}(c_{j_1}, \ldots, c_{j_n}).$$
\end{Lemma1}

\begin{proof}
Let $\mathcal{M} \models \mathbf{W}_1$. Let $\mathbf{t}(x_1, \ldots, x_n)$ be an $\mathcal{L}^*$ Skolem function and, let $i_1 < \ldots < i_n$ and $j_1 < \ldots < j_n$ be in $\mathbb{Z}$. Assume that for all $k \in \mathbb{Z}$, $\mathcal{M} \models \mathbf{t}(c_{i_1}, \ldots, c_{i_n}) < c_k$. Let $k \in \mathbb{Z}$ be such that $k < i_1$. By indiscernibility (Definition \ref{Df:TheoryW1}(vii)):
$$\mathcal{M} \models \forall x(\mathcal{C}(x) \Rightarrow \mathbf{t}(c_{i_1}, \ldots, c_{i_n}) < x).$$
Therefore by Definition \ref{Df:TheoryW1}(viii):
$$\mathcal{M} \models \mathbf{t}(c_{i_1}, \ldots, c_{i_n})= \mathbf{t}(c_{j_1}, \ldots, c_{j_n}).$$ 
\Square
\end{proof}

\begin{Theorems1} \label{Th:AutomorphismMovingNoNaturalDown}
If $\mathcal{M}= \langle M^\mathcal{M}, \in^\mathcal{M} \rangle$ is an $\omega$-model of $\mathrm{Mac}$ then there exists an $\mathcal{L}$-structure $\mathcal{N} \equiv \mathcal{M}$ admitting an automorphism $j: \mathcal{N} \longrightarrow \mathcal{N}$ such that
\begin{itemize}
\item[(i)] $\mathcal{N} \models j(n) \geq n$ for all $n \in (\omega^\mathcal{N})^*$,
\item[(ii)] there exists an $n \in (\omega^\mathcal{N})^*$ with $\mathcal{N} \models j(n) > n$. 
\end{itemize}  
\end{Theorems1}

\begin{proof}
Let $\mathcal{M}= \langle M^\mathcal{M}, \in^\mathcal{M} \rangle$ be an $\omega$-model of $\mathrm{Mac}$. By Lemma \ref{Th:ConsistencyOfW1} there exists an $\mathcal{L}^{\mathbf{ind}}$-structure $\mathcal{Q} \models \mathbf{W}_1$ such that $\mathcal{Q} \mid_{\mathcal{L}} \equiv \mathcal{M}$. Let $C= \{ c_i^\mathcal{Q} \mid i \in \mathbb{Z} \}$. Let $\mathcal{N}= \mathcal{H}_{\mathcal{L}}^{\mathcal{Q}}(C)$. Therefore $\mathcal{N} \equiv \mathcal{M}$. Define $j: \mathcal{N} \longrightarrow \mathcal{N}$ by 
$$j(c_i^\mathcal{Q})= c_{i+1}^\mathcal{Q} \textrm{ for all } i \in \mathbb{Z},$$
$$\textrm{for all } \mathcal{L}^* \textrm{ Skolem functions } \mathbf{t}(x_1, \ldots, x_n) \textrm{ and for all } i_1 < \ldots < i_n \textrm{ in } \mathbb{Z},$$
$$j(\mathbf{t}(c_{i_1}^\mathcal{Q}, \ldots, c_{i_n}^\mathcal{Q}))= \mathbf{t}(j(c_{i_1}^\mathcal{Q}), \ldots, j(c_{i_n}^\mathcal{Q})).$$
Now, $j$ is an automorphism of the $\mathcal{L}$-structure $\mathcal{N}$. If $\mathbf{t}(x_1, \ldots x_n)$ is an $\mathcal{L}^*$ Skolem function and $i_1 < \ldots < i_n$ is in $\mathbb{Z}$ such that $\mathcal{Q} \models \mathbf{t}(c_{i_1}, \ldots, c_{i_n}) \in \omega$ then by Lemma \ref{Th:IndiscerniblesCofinalInSkolemHull} either
\begin{itemize}
\item[(i)] there is $k \in \mathbb{Z}$ such that 
$$\mathcal{Q} \models c_k < \mathbf{t}(c_{i_1}, \ldots, c_{i_n}) \leq c_{k+1}$$
(Note that $k$ might equal $i_j$ for some $1 \leq j \leq n$),  
\item[(ii)] or, for all $k \in \mathbb{Z}$,
$$\mathcal{Q} \models \mathbf{t}(c_{i_1}, \ldots, c_{i_n}) < c_k.$$  
\end{itemize}
If (ii) holds then Lemma \ref{Th:SkolemFunctionsFixedBelow} implies that
$$j(\mathbf{t}(c_{i_1}^\mathcal{Q}, \ldots, c_{i_n}^\mathcal{Q}))= \mathbf{t}(c_{i_1+1}^\mathcal{Q}, \ldots, c_{i_n+1}^\mathcal{Q})= \mathbf{t}(c_{i_1}^\mathcal{Q}, \ldots, c_{i_n}^\mathcal{Q}).$$
If (i) holds then indiscernibility (Definition \ref{Df:TheoryW1}(vii)) implies that
$$\mathcal{Q} \models c_k < \mathbf{t}(c_{i_1}, \ldots, c_{i_n}) \leq c_{k+1} \textrm{ if and only if } \mathcal{Q} \models c_{k+1} < \mathbf{t}(c_{i_1+1}, \ldots, c_{i_n+1}) \leq c_{k+2}.$$
$$\textrm{Therefore } j(\mathbf{t}(c_{i_1}^\mathcal{Q}, \ldots, c_{i_n}^\mathcal{Q}))= \mathbf{t}(c_{i_1+1}^\mathcal{Q}, \ldots, c_{i_n+1}^\mathcal{Q}) >^\mathcal{Q} \mathbf{t}(c_{i_1}^\mathcal{Q}, \ldots, c_{i_n}^\mathcal{Q}).$$
This shows that $\mathcal{N} \models j(n) \geq n$ for all $n \in (\omega^\mathcal{N})^*$. And, $\mathcal{N} \models j(c_0^\mathcal{Q}) > c_0^\mathcal{Q}$. 
\Square
\end{proof}

Theorem \ref{Th:AutomorphismMovingNoNaturalDown} can be extended to show that the existence of a transitive model of Mac Lane set theory implies the existence of an elementarily equivalent model admitting a non-trivial automorphism that does not move any ordinal down. We begin by defining an extension of the $\mathcal{L}^{\mathbf{ind}}$-theory $\mathbf{W}_1$.

\begin{Definitions1} \label{Df:TheoryW2}
We define the $\mathcal{L}^{\mathbf{ind}}$-theory $\mathbf{W}_2 \supseteq \mathrm{Mac} \cup \mathbf{Sk}(\mathcal{L})$ by adding the axioms
\begin{itemize}
\item[(i)] $\forall x (\mathcal{C}(x) \Rightarrow x \in \omega)$
\item[(ii)] $(\forall x \in \omega)(\exists y \in \omega)(y > x \land \mathcal{C}(y))$,
\item[(iii)] $\mathcal{C}(c_i)$ for all $i \in \mathbb{Z}$,
\item[(iv)] $c_i < c_j$ for all $i, j \in \mathbb{Z}$ with $i < j$,
\item[(v)] for all $\mathcal{L}^*$-formulae $\phi(x_1, \ldots, x_n)$,
$$\forall x_1 \ldots \forall x_n \forall y_1 \ldots \forall y_n \left( 
\begin{array}{c}
x_1 < \ldots < x_n \land y_1 < \ldots < y_n \land \bigwedge_{1 \leq i \leq n} (\mathcal{C}(x_i) \land \mathcal{C}(y_i))\\
\Rightarrow (\phi(x_1, \ldots, x_n) \iff \phi(y_1, \ldots, y_n))
\end{array} \right),$$
\item[(vi)] for all $\mathcal{L}^*$ Skolem functions $\mathbf{t}(x_1, \ldots, x_n)$ and for all $i_1 < \ldots < i_n$ and $j_1 < \ldots < j_n$ in $\mathbb{Z}$,
$$\forall x (\mathcal{C}(x) \Rightarrow \mathbf{t}(c_{i_1}, \ldots, c_{i_n}) < x) \Rightarrow (\mathbf{t}(c_{i_1}, \ldots, c_{i_n})= \mathbf{t}(c_{j_1}, \ldots, c_{j_n})),$$
\item[(vii)] for all $\mathcal{L}^*$ Skolem functions $\mathbf{t}(x_1, \ldots, x_n)$ and for all $i_1 < \ldots < i_n$ in $\mathbb{Z}$,
$$(\mathbf{t}(c_{i_1}, \ldots, c_{i_n}) \textrm{ is an ordinal}) \Rightarrow (\mathbf{t}(c_{i_1}, \ldots, c_{i_n}) \leq \mathbf{t}(c_{i_n+1}, \ldots, c_{i_n+1})).$$   
\end{itemize}
\end{Definitions1}
            
\noindent It should be noted that (vii) of Definition \ref{Df:TheoryW2} is an axiom scheme. This scheme will ensure that the ordinals in the Skolem hull of a model of $\mathbf{W}_2$ that includes the order indiscernible $c_i$s will not be moved down by the automorphism of this Skolem hull generated by sending $c_i$ to $c_{i+1}$ for all $i \in \mathbb{Z}$.\\
\\
\indent We are now in a position to prove an analogue of Lemma \ref{Th:ConsistencyOfW1}.

\begin{Lemma1} \label{Th:ConsistencyOfW2}
If $\langle M, \in \rangle$ is a transitive model of $\mathrm{Mac}$ then $\mathbf{W}_2 \cup \mathbf{Th}_\mathcal{L}(\langle M, \in \rangle)$ is consistent.
\end{Lemma1}

\begin{proof}
Let $\mathcal{M}= \langle M, \in \rangle$ be a transitive model of $\mathrm{Mac}$. Therefore $\langle M, \in \rangle$ is an $\omega$-model and the $\omega$ of $\langle M, \in \rangle$ coincides with the $\omega$ of the ambient model. Again, we will use the Compactness Theorem to show that $\mathbf{W}_2 \cup \mathbf{Th}_\mathcal{L}(\langle M, \in \rangle)$ is consistent. Let $\Delta \subseteq \mathbf{W}_2$ be finite. Using the same arguments used in the proof of Lemma \ref{Th:ConsistencyOfW1} we can find $C \subseteq \omega$ and expand $\mathcal{M}$ to a structure $\mathcal{M}^\prime$ in which $\mathcal{C}$ is interpreted using $C$ and
\begin{itemize}
\item[(i)] $\mathcal{M}^\prime \models \mathbf{Sk}(\mathcal{L})$,  
\item[(ii)] $\mathcal{M}^\prime \mid_{\mathcal{L}}= \mathcal{M}$, 
\item[(iii)] $\mathcal{M}^\prime$ satisfies (i) and (ii) of Definition \ref{Df:TheoryW2},
\item[(iv)] $\mathcal{M}^\prime$ satisfies all instances of Definition \ref{Df:TheoryW2}(v) that appear in $\Delta$. 
\end{itemize}
We can also use the arguments utilised in the proof of Lemma \ref{Th:ConsistencyOfW1} to find an infinite $D \subseteq C$ such that if $\mathbf{t}(x_1, \ldots, x_k)$ is a Skolem function mentioned in $\Delta$, and $a_1 < \ldots < a_k$ and $b_1 < \dots < b_k$ are in $D$ then
$$\textrm{if } \mathbf{t}^{\mathcal{M}^\prime}(a_1, \ldots, a_k) \in \omega \textrm{ and } \mathbf{t}^{\mathcal{M}^\prime}(a_1, \ldots, a_k) < \min C$$ 
$$\textrm{then } \mathbf{t}^{\mathcal{M}^\prime}(a_1, \ldots, a_k)= \mathbf{t}^{\mathcal{M}^\prime}(b_1, \ldots, b_k).$$
Now, suppose that $\Delta$ only mentions the constant symbols $c_{-n}, \ldots, c_n$. Let $\mathbf{s}_0, \ldots, \mathbf{s}_{m-1}$ and $\eta_0, \ldots, \eta_{m-1}$ be such that
\begin{itemize}
\item[(i)] $\eta_i \in \mathbb{Z}^{<\omega}$ for all $i \in m$, 
\item[(ii)] $\eta_i(p) < \eta_i(q)$ for all $i \in m$ and $p < q < |\eta_i|$,
\item[(iii)] $\mathbf{s}_i$ is an $\mathcal{L}^*$ Skolem function of arity $|\eta_i|$ for all $i \in m$,
\item[(iv)] if the axiom
$$(\mathbf{t}(c_{i_1}, \ldots, c_{i_l}) \textrm{ is an ordinal}) \Rightarrow (\mathbf{t}(c_{i_1}, \ldots, c_{i_l}) \leq \mathbf{t}(c_{i_1+1}, \ldots, c_{i_l+1})) \textrm{ appears in } \Delta$$
then there is an $j \in m$ such that $\mathbf{t}= \mathbf{s}_j$ and $\eta_j= \langle i_1, \ldots, i_l \rangle$.
\end{itemize}
Note that if $i_1 < \ldots < i_l$ are in $\mathbb{Z}$ then $i_1 < i_1+1 \leq i_2 < i_2+1 \leq \ldots i_{l-1}+1 \leq i_l$. For each $0 \leq j < m$, let $\psi_j$ be the $\mathcal{L}^*$-formula obtained by replacing $c_{i_1}, c_{i_1+1}, \ldots, c_{i_l}, c_{i_l+1}$ by free variables $x_1, \ldots, x_{u_j}$ in the $\mathcal{L}^{\mathbf{ind}}$-sentence
$$\mathbf{s}_j(c_{i_1}, \ldots, c_{i_l}) \leq \mathbf{s}_j(c_{i_1+1}, \ldots, c_{i_l+1}) \textrm{ where } \eta_j= \langle i_1, \ldots, i_l \rangle.$$
Note that for all $0 \leq j < m$, the $\mathcal{L}^*$-formula $\psi_j$ has arity $u_j$ for $|\eta_j|+ 1 \leq u_j \leq 2\cdot|\eta_j|$. For all $j \in m$ and $k \in |\eta_j|$ we define $\zeta_k^j$ by
$$\zeta_0^j= 1 \textrm{ for all } j \in m,$$
$$\zeta_{k+1}^j= \left\{ \begin{array}{ll}
\zeta_k^j+1 & \textrm{if } \eta_j(k+1) = \eta_j(k) + 1\\
\zeta_k^j+2 & \textrm{otherwise}
\end{array}\right)$$
This ensures that if $j \in m$ and $x_1, \ldots, x_{u_j}$ are the free variables appearing in $\psi_j$ then
$$\psi_j(x_1, \ldots, x_{u_j}) \textrm{ is } \mathbf{s}_j(x_{\zeta_0^j}, \ldots, x_{\zeta_{|\eta_j|-1}^j}) \leq \mathbf{s}_j(x_{\zeta_0^j+1}, \ldots, x_{\zeta_{|\eta_j|-1}^j+1}).$$  
Now, consider the function $K_0: [D]^{u_0} \longrightarrow 3$ defined by
$$K_0(\{x_1, \ldots, x_{u_0}\})= \left\{ \begin{array}{ll}
0 & \textrm{if } \mathcal{M}^{\prime} \models \left( \begin{array}{c}
(\mathbf{s}_0(x_{\zeta_0^0}, \ldots, x_{\zeta_{|\eta_0|-1}^0}) \textrm{ is an ordinal})\\
\land (\neg \psi_0(x_1, \ldots, x_{u_0}))
\end{array} \right) \\
1 & \textrm{if } \mathcal{M}^{\prime} \models \left( \begin{array}{c}
(\mathbf{s}_0(x_{\zeta_0^0}, \ldots, x_{\zeta_{|\eta_0|-1}^0}) \textrm{ is an ordinal})\\
\land(\psi_0(x_1, \ldots, x_{u_0}))
\end{array}\right)\\
2 & \textrm{otherwise} 
\end{array}\right)$$
$$\textrm{where } x_1 < \ldots < x_{u_0}.$$
By Ramsey's Theorem there is an infinite $A_0 \subseteq D$ such that $K_0``[A_0]^{u_0}= \{p\}$ for some $p \in 3$. Now, for each $0 \leq j < m-1$ define $K_{j+1}: [A_j]^{u_{j+1}} \longrightarrow 3$ by
$$K_{j+1}(\{x_1, \ldots, x_{u_{j+1}}\})= \left\{ \begin{array}{ll}
0 & \textrm{if } \mathcal{M}^{\prime} \models \left( \begin{array}{c}
(\mathbf{s}_{j+1}(x_{\zeta_0^{j+1}}, \ldots, x_{\zeta_{|\eta_{j+1}|-1}^{j+1}}) \textrm{ is an ordinal})\\
\land(\neg \psi_{j+1}(x_1, \ldots, x_{u_{j+1}}))
\end{array} \right)\\
1 & \textrm{if } \mathcal{M}^{\prime} \models \left( \begin{array}{c}
(\mathbf{s}_{j+1}(x_{\zeta_0^{j+1}}, \ldots, x_{\zeta_{|\eta_{j+1}|-1}^{j+1}}) \textrm{ is an ordinal})\\
 \land(\psi_{j+1}(x_1, \ldots, x_{u_{j+1}}))
\end{array} \right) \\
2 & \textrm{otherwise} 
\end{array}\right)$$
$$\textrm{where } x_1 < \ldots < x_{u_{j+1}}.$$
And let $A_{j+1} \subseteq A_j$ be the infinite set guaranteed by Ramsey's Theorem such that $K_{j+1}``[A_{j+1}]^{u_{j+1}}= \{p\}$ for some $p \in 3$.\\
\textbf{Claim: } For all $0 \leq j < m$, $K_j``[A_{m-1}]^{u_j} \neq \{ 0 \}$.\\
Suppose that $j \in m$ with $K_j``[A_{m-1}]^{u_j} = \{ 0 \}$. Let $h: \omega \longrightarrow A_{m-1}$ enumerate the elements of $A_{m-1}$ in order.
$$\textrm{Let } B= \{ \mathbf{s}_j^{\mathcal{M}^{\prime}}(h(\zeta_0^j+k), \ldots, h(\zeta_{|\eta_j|-1}^j+k)) \mid k \in \omega \}.$$
Now, $B$ is a set of ordinals. Let $\alpha$ be the least element of $B$. Therefore there is a $k \in \omega$ such that  $\alpha= \mathbf{s}_j^{\mathcal{M}^{\prime}}(h(\zeta_0^j+k), \ldots, h(\zeta_{|\eta_j|-1}^j+k))$. But now,
$$\mathcal{M}^{\prime} \models \mathbf{s}_j(h(\zeta_0^j+k), \ldots, h(\zeta_{|\eta_j|-1}^j+k)) > \mathbf{s}_j(h(\zeta_0^j+(k+1)), \ldots, h(\zeta_{|\eta_j|-1}^j+(k+1)))$$
$$\textrm{and } \mathbf{s}_j^{\mathcal{M}^{\prime}}(h(\zeta_0^j+(k+1)), \ldots, h(\zeta_{|\eta_j|-1}^j+(k+1))) \in B$$
which contradicts the fact that $\mathbf{s}_j^{\mathcal{M}^{\prime\prime}}(h(\zeta_0^j+k), \ldots, h(\zeta_{|\eta_j|-1}^j+k))$ is least. This proves the claim.\\
Using the first $2\cdot n +1$ elements of $A_{m-1}$ to interpret the constants $c_{-n}, \ldots, c_n$ we can expand $\mathcal{M}^{\prime}$ to a structure $\mathcal{M}^{\prime\prime}$ such that
$$\mathcal{M}^{\prime\prime} \models \Delta \cup \mathbf{Th}_{\mathcal{L}}(\langle M, \in \rangle) \cup \mathbf{Sk}(\mathcal{L}).$$
Therefore, by compactness, $\mathbf{W}_2 \cup \mathbf{Th}_{\mathcal{L}}(\langle M, \in \rangle)$ is consistent. 
\Square
\end{proof}

Since the theory $\mathbf{W}_2$ extends $\mathbf{W}_1$, Lemmas \ref{Th:IndiscerniblesCofinalInSkolemHull} and \ref{Th:SkolemFunctionsFixedBelow} apply to models of $\mathbf{W}_2$. This allows us to prove an analogue of Theorem \ref{Th:AutomorphismMovingNoNaturalDown} showing that if there exists a transitive model of $\mathrm{Mac}$ then there exists an elementarily equivalent model admitting a non-trivial automorphism that does not move any ordinal down.

\begin{Theorems1} \label{Th:TransitiveModelsYieldAutomorphismMovingNoOrdinalDown}
If $\langle M, \in \rangle$ is a transitive model of $\mathrm{Mac}$ then there is an $\mathcal{L}$-structure $\mathcal{N} \equiv \langle M, \in \rangle$ admitting an automorphism $j:\mathcal{N} \longrightarrow \mathcal{N}$ such that 
\begin{itemize}
\item[(i)] $\mathcal{N} \models j(\alpha) \geq \alpha$ for all $\alpha \in \mathrm{Ord}^\mathcal{N}$, 
\item[(ii)] there exists an $n \in (\omega^\mathcal{N})^*$ such that $\mathcal{N} \models j(n) > n$.
\end{itemize}
\end{Theorems1}

\begin{proof}
Let $\langle M, \in \rangle$ be a transitive model of $\mathrm{Mac}$. By Lemma \ref{Th:ConsistencyOfW2} there exists an $\mathcal{L}^{\mathbf{ind}}$-structure $\mathcal{Q} \models \mathbf{W}_2$ such that $\mathcal{Q} \mid_{\mathcal{L}} \equiv \langle M, \in \rangle$. Let $C= \{ c_i^\mathcal{Q} \mid i \in \mathbb{Z} \}$. Let $\mathcal{N}= \mathcal{H}_{\mathcal{L}}^{\mathcal{Q}}(C)$. Therefore $\mathcal{N} \equiv \langle M, \in \rangle$. Define $j: \mathcal{N} \longrightarrow \mathcal{N}$ by 
$$j(c_i^\mathcal{Q})= c_{i+1}^\mathcal{Q} \textrm{ for all } i \in \mathbb{Z},$$
$$\textrm{for all } \mathcal{L}^* \textrm{ Skolem functions } \mathbf{t}(x_1, \ldots, x_n) \textrm{ and for all } i_1 < \ldots < i_n \textrm{ in } \mathbb{Z},$$
$$j(\mathbf{t}(c_{i_1}^\mathcal{Q}, \ldots, c_{i_n}^\mathcal{Q}))= \mathbf{t}(j(c_{i_1}^\mathcal{Q}), \ldots, j(c_{i_n}^\mathcal{Q})).$$
Now, $j$ is an automorphism of the $\mathcal{L}$-structure $\mathcal{N}$. The same arguments used in the proof of Theorem \ref{Th:AutomorphismMovingNoNaturalDown} reveal that if $\mathbf{t}(x_1, \ldots, x_n)$ is an $\mathcal{L}^*$ Skolem function and $i_1 < \ldots < i_n$ are in $\mathbb{Z}$ such that $\mathcal{N} \models \mathbf{t}(c_{i_1}^\mathcal{Q}, \ldots, c_{i_n}^\mathcal{Q}) \in \omega$ then 
$$\mathcal{N} \models j(\mathbf{t}(c_{i_1}^\mathcal{Q}, \ldots, c_{i_n}^\mathcal{Q}))= \mathbf{t}(c_{i_1}^\mathcal{Q}, \ldots, c_{i_n}^\mathcal{Q}) \textrm{ or } \mathcal{N} \models j(\mathbf{t}(c_{i_1}^\mathcal{Q}, \ldots, c_{i_n}^\mathcal{Q})) > \mathbf{t}(c_{i_1}^\mathcal{Q}, \ldots, c_{i_n}^\mathcal{Q}).$$
In particular,
$$\mathcal{N} \models j(c_0^\mathcal{Q}) > c_0^\mathcal{Q}.$$
If $\mathbf{t}(x_1, \ldots, x_n)$ is an $\mathcal{L}^*$ Skolem function and $i_1 < \ldots < i_n$ are in $\mathbb{Z}$ and
$$\mathcal{Q} \models \mathbf{t}(c_{i_1}, \ldots, c_{i_n}) \textrm{ is an ordinal}$$ 
then the scheme (vii) of Definition \ref{Df:TheoryW2} implies that
$$j(\mathbf{t}^\mathcal{Q}(c_{i_1}^\mathcal{Q}, \ldots, c_{i_n}^\mathcal{Q}))= \mathbf{t}^\mathcal{Q}(c_{i_1+1}^\mathcal{Q}, \ldots, c_{i_n+1}^\mathcal{Q}) \geq \mathbf{t}^\mathcal{Q}(c_{i_1}^\mathcal{Q}, \ldots, c_{i_n}^\mathcal{Q}).$$
Therefore for all $\alpha \in \mathrm{Ord}^\mathcal{N}$,
$$\mathcal{N} \models j(\alpha) \geq \alpha.$$ 
\Square
\end{proof}

\noindent This theorem allows us to prove that every complete consistent extension of $\mathrm{ZFC}$ has a model admitting a non-trivial automorphism that does not move any ordinal down. Using the Reflection Principle and the Compactness Theorem we show that every complete consistent extension of $\mathrm{ZFC}$ has a model $\mathcal{M}$ with a point $\lambda \in M^\mathcal{M}$, which $\mathcal{M}$ believes is a limit ordinal, such that externally the structure $\langle V_\lambda^\mathcal{M}, \in^\mathcal{M} \rangle$ is elementarily equivalent to $\mathcal{M}$. Since $\mathcal{M}$ believes that $V_\lambda$ is a transitive model of $\mathrm{Mac}$, we can apply Lemma \ref{Th:ConsistencyOfW2} inside $\mathcal{M}$ to obtain a model admitting the desired automorphism.

\begin{Theorems1} \label{Th:ZFCwithAutThatDoesNotMoveAnyOrdinalDown}
If $T \supseteq \mathrm{ZFC}$ is complete and consistent $\mathcal{L}$-theory then there is an $\mathcal{L}$-structure $\mathcal{N} \models T$ admitting an automorphism $j: \mathcal{N} \longrightarrow \mathcal{N}$ such that 
\begin{itemize}
\item[(i)] $\mathcal{N} \models j(\alpha) \geq \alpha$ for all $\alpha \in \mathrm{Ord}^\mathcal{N}$, 
\item[(ii)] there exists an ordinal $\alpha \in \mathrm{Ord}^\mathcal{N}$ such that $\mathcal{N} \models j(\alpha) > \alpha$. 
\end{itemize}
\end{Theorems1}

\begin{proof}
Let $T \supseteq \mathrm{ZFC}$ be complete and consistent $\mathcal{L}$-theory. Let $\bar{\lambda}$ be a new constant symbol. Define $T^\star \supseteq T$ to be the $\mathcal{L}_{\bar{\lambda}}$-theory that includes the axioms:
\begin{itemize}
\item[(I)] $\bar{\lambda}$ is a limit ordinal, 
\item[(II)] for all $\mathcal{L}$-sentences, $\sigma$,
$$(\langle V_{\bar{\lambda}}, \in \rangle \models \sigma) \iff \sigma.$$ 
\end{itemize} 
Let $\mathcal{M}$ be an $\mathcal{L}$-structure such that $\mathcal{M} \models T$. Let $\Delta \subseteq T^\star$ be finite. Using the Reflection Principle we can expand $\mathcal{M}$ to an $\mathcal{L}_{\bar{\lambda}}$-structure $\mathcal{M}^\prime$ such that $\mathcal{M}^\prime \models \Delta$. Therefore, by the Compactness Theorem, there is an $\mathcal{L}_{\bar{\lambda}}$-structure $\mathcal{Q}$ such that
$$\mathcal{Q} \models T^\star.$$
We work inside $\mathcal{Q}$. Since $\bar{\lambda}$ is a limit ordinal, $V_{\bar{\lambda}}$ is a transitive model of $\mathrm{Mac}$. Therefore Theorem \ref{Th:TransitiveModelsYieldAutomorphismMovingNoOrdinalDown} implies that there exists an $\mathcal{L}$-structure $\mathcal{N} \equiv \langle V_{\bar{\lambda}}, \in \rangle$ admitting an automorphism $j: \mathcal{N} \longrightarrow \mathcal{N}$ such that
\begin{itemize}
\item[(i)] $\mathcal{N} \models j(\alpha) \geq \alpha$ for all $\alpha \in \mathrm{Ord}^\mathcal{N}$, 
\item[(ii)] there exists an ordinal $\alpha \in \mathrm{Ord}^\mathcal{N}$ such that $\mathcal{N} \models j(\alpha) > \alpha$. 
\end{itemize}
Axiom scheme (II) above ensures that $\mathcal{N} \models T$.   
\Square
\end{proof}

\section[Automorphisms from non-standard models]{Automorphisms from non-standard models} \label{Sec:AutomorphismsFromNonStandardModels}

In Theorem \ref{Th:AutomorphismMovingNoNaturalDown} and Theorem \ref{Th:TransitiveModelsYieldAutomorphismMovingNoOrdinalDown} we built models that admit automorphisms which are well-behaved in the sense that they do not move any point down along an initial segment of the ordinals. In this section we will use non-standard models of set theory with standard parts to construct models of set theory equipped with automorphisms which do not move any points down along an initial segment of the ordinals but do move points down above this initial segment. This allows us to answer Question \ref{Q:HolmesQuestion}. \\
\\
\indent We begin by showing that models of set theory admitting an automorphism which does not move any natural number down and does move an ordinal above the natural numbers down can be built from a non-standard $\omega$-model of set theory. Throughout this section $\hat{f}$ will be a new unary function symbol and $\bar{\beta}$ will be a new constant symbol. This allows us to define the language $\mathcal{L}^{\mathbf{ind}}_{\hat{f}, \bar{\beta}}$ that extends $\mathcal{L}^*_{\hat{f}, \bar{\beta}}$ in the same way that $\mathcal{L}^{\mathbf{ind}}$ extends $\mathcal{L}^*$.

\begin{Definitions1}
We define $\mathcal{L}^{\mathbf{ind}}_{\hat{f}, \bar{\beta}} \supseteq \mathcal{L}^*_{\hat{f}, \bar{\beta}}$ by adding
\begin{itemize}
\item[(i)] a unary predicate $\mathcal{C}$,
\item[(ii)] constant symbols $c_i$ for each $i \in \mathbb{Z}$.
\end{itemize}
\end{Definitions1}

\noindent Equipped with this language we are able to define a theory extending $\mathbf{W}_1$ that yields models admitting the desired automorphism.    

\begin{Definitions1} \label{Df:TheoryW3}
We define the $\mathcal{L}^{\mathbf{ind}}_{\hat{f}, \bar{\beta}}$-theory $\mathbf{W}_3 \supseteq \mathrm{Mac} \cup \mathbf{Sk}(\mathcal{L}_{\hat{f}, \bar{\beta}})$ by adding the axioms
\begin{itemize}
\item[(i)] $\forall x (\mathcal{C}(x) \Rightarrow x \in \omega)$
\item[(ii)] $(\forall x \in \omega)(\exists y \in \omega)(y > x \land \mathcal{C}(y))$,
\item[(iii)] $\mathcal{C}(c_i)$ for all $i \in \mathbb{Z}$,
\item[(iv)] $c_i < c_j$ for all $i, j \in \mathbb{Z}$ with $i < j$,
\item[(v)] $\bar{\beta}$ is an ordinal,
\item[(vi)] $(\forall x \in \omega)(\hat{f}(x) \in \bar{\beta})$,
\item[(vii)] $(\forall x \in \omega)(\forall y \in \omega)(x < y \Rightarrow \hat{f}(y) < \hat{f}(x))$,
\item[(viii)] for all $\mathcal{L}^*_{\hat{f}, \bar{\beta}}$-formulae $\phi(x_1, \ldots, x_n)$,
$$\forall x_1 \ldots \forall x_n \forall y_1 \ldots \forall y_n \left( 
\begin{array}{c}
x_1 < \ldots < x_n \land y_1 < \ldots < y_n \land \bigwedge_{1 \leq i \leq n} (\mathcal{C}(x_i) \land \mathcal{C}(y_i))\\
\Rightarrow (\phi(x_1, \ldots, x_n) \iff \phi(y_1, \ldots, y_n))
\end{array} \right),$$
\item[(ix)] for all $\mathcal{L}^*_{\hat{f}, \bar{\beta}}$ Skolem functions $\mathbf{t}(x_1, \ldots, x_n)$ and for all $i_1 < \ldots < i_n$ and \mbox{$j_1 < \ldots < j_n$} in $\mathbb{Z}$,
$$\forall x (\mathcal{C}(x) \Rightarrow \mathbf{t}(c_{i_1}, \ldots, c_{i_n}) < x) \Rightarrow (\mathbf{t}(c_{i_1}, \ldots, c_{i_n})= \mathbf{t}(c_{j_1}, \ldots, c_{j_n})).$$   
\end{itemize}
\end{Definitions1}

As we saw with the theory $\mathbf{W}_1$, axioms (i), (ii) and (iii) and axiom scheme (vii) of Definition \ref{Df:TheoryW3} ensure that $\mathcal{C}$ defines a class of order indiscernibles that is cofinal in the natural numbers and the $c_i$s form a subclass of $\mathcal{C}$. Axiom scheme (ix) ensures that every Skolem term which is interpreted by a natural number less than every element in $\mathcal{C}$ is fixed when the sequence of $c_i$s appearing in this term is replaced by an order equivalent sequence of $c_i$s. Axioms (v), (vi) and (vii) ensure that $\hat{f}$ is an order reversing function from the natural numbers into the ordinal $\bar{\beta}$. As with $\mathbf{W}_1$, the Skolem hull of a model of $\mathbf{W}_3$ that contains the $c_i$s will admit an automorphism that does not move any natural number down. The behaviour of the unary function symbol $\hat{f}$ will ensure that any automorphism of a substructure of a model of $\mathbf{W}_3$ moves the point interpreting $\hat{f}(c_0)$ in the opposite direction to the point interpreting $c_0$.\\      
\\    
\indent We are now in a position to show that the existence of an $\omega$-model of $\mathrm{Mac}$ that is non-standard below an ordinal $\bar{\beta}$ implies the existence of an elementarily equivalent model of $\mathbf{W}_3$.

\begin{Lemma1} \label{Th:ConsistencyOfW3}
If $\mathcal{M}= \langle M^\mathcal{M}, \in^\mathcal{M}, \bar{\beta}^\mathcal{M} \rangle$ is an $\omega$-model of $\mathrm{Mac}$ such that\\ 
$\mathcal{M} \models \bar{\beta} \textrm{ is an ordinal}$, and $\mathcal{M}$ is non-standard below $\bar{\beta}^\mathcal{M}$ then $\mathbf{W}_3 \cup \mathbf{Th}_{\mathcal{L}_{\bar{\beta}}}(\mathcal{M})$ is consistent.
\end{Lemma1}

\begin{proof}
Let $\mathcal{M}= \langle M^\mathcal{M}, \in^\mathcal{M}, \bar{\beta}^\mathcal{M} \rangle$ be an $\omega$-model of $\mathrm{Mac}$ such that 
$$\mathcal{M} \models \bar{\beta} \textrm{ is an ordinal}$$ 
and $\mathcal{M}$ is non-standard below $\bar{\beta}^\mathcal{M}$. Let $g: \omega \longrightarrow (\omega^\mathcal{M})^*$ be the isomorphism guaranteed by the fact that $\mathcal{M}$ is an $\omega$-model. Let $f: \omega \longrightarrow (\bar{\beta}^\mathcal{M})^*$ be such for all $n, m \in \omega$, 
$$n < m \textrm{ if and only if } f(m) <^\mathcal{M} f(n).$$
The existence of such a function is guaranteed by the fact that $\mathcal{M}$ is non-standard below $\bar{\beta}^\mathcal{M}$. As with the proof of Lemma \ref{Th:ConsistencyOfW1} and Lemma \ref{Th:ConsistencyOfW2} we will use compactness to show that $\mathbf{W}_3 \cup \mathbf{Th}_{\mathcal{L}_{\bar{\beta}}}$ is consistent. Let $\Delta \subseteq \mathbf{W}_3$ be finite. Using $f$ to interpret $\hat{f}$ we can expand $\mathcal{M}$ to an $\mathcal{L}_{\hat{f}, \bar{\beta}}$-structure $\mathcal{M}^\prime$ that satisfies (v), (vi) and (vii) of Definition \ref{Df:TheoryW3}. Using the same arguments we used in the proof of Lemma \ref{Th:ConsistencyOfW1} we can find $C \subseteq (\omega^\mathcal{M})^*$ and expand $\mathcal{M}^\prime$ to a structure $\mathcal{M}^{\prime\prime}$ in which $\mathcal{C}$ is interpreted using $C$ and
\begin{itemize}
\item[(i)] $\mathcal{M}^{\prime\prime} \models \mathbf{Sk}(\mathcal{L}_{\hat{f}, \bar{\beta}})$, 
\item[(ii)] $\mathcal{M}^{\prime\prime} \mid_{\mathcal{L}_{\hat{f}, \bar{\beta}}}= \mathcal{M}^\prime$, 
\item[(iii)] $\mathcal{M}^{\prime\prime}$ satisfies (i) and (ii) of Definition \ref{Df:TheoryW3},
\item[(iv)] $\mathcal{M}^{\prime\prime}$ satisfies all instances of \ref{Df:TheoryW3}(viii) appearing in $\Delta$. 
\end{itemize}
Without loss of generality we can assume that $\Delta$ only mentions the constant symbols $c_{-n}, \ldots, c_0, \ldots c_n$ and $\bar{\beta}$. Again, using the same arguments we used in the proof Lemma \ref{Th:ConsistencyOfW1} we can find an infinite $D \subseteq C$ such that by interpreting $c_{-n}, \ldots, c_0, \ldots, c_n$ using the first $2\cdot n +1$ of $D$ we can expand $\mathcal{M}^{\prime\prime}$ to a structure $\mathcal{M}^{\prime\prime\prime}$ that satisfies all instance of (iii) and (ix) of Definition \ref{Df:TheoryW3} which appear in $\Delta$. Therefore
$$\mathcal{M}^{\prime\prime\prime} \models \Delta \cup \mathbf{Th}_{\mathcal{L}_{\bar{\beta}}}(\mathcal{M}) \cup \mathbf{Sk}(\mathcal{L}_{\hat{f}, \bar{\beta}}).$$
Therefore, by compactness, $\mathbf{W}_3 \cup \mathbf{Th}_{\mathcal{L}_{\bar{\beta}}}(\mathcal{M})$ is consistent.  
\Square
\end{proof}

Axiom schemes (viii) and (ix) of Definition \ref{Df:TheoryW3} allow us to prove analogues of Lemmas \ref{Th:IndiscerniblesCofinalInSkolemHull} and \ref{Th:SkolemFunctionsFixedBelow}.

\begin{Lemma1} \label{Th:IndiscerniblesCofinalInSkolemHull2}
(K\"{o}rner \cite{kor94}) Let $\mathcal{M} \models \mathbf{W}_3$. For all $\mathcal{L}^*_{\hat{f}, \bar{\beta}}$ Skolem functions $\mathbf{t}(x_1, \ldots, x_n)$ and for all $i_1 < \ldots < i_n$ in $\mathbb{Z}$, if $\mathcal{M} \models \mathbf{t}(c_{i_1}, \ldots, c_{i_n}) \in \omega$ then there exists a $k \in \mathbb{Z}$ such that
$$\mathcal{M} \models \mathbf{t}(c_{i_1}, \ldots, c_{i_n}) < c_k.$$  
\end{Lemma1}

\begin{proof}
Using the same arguments as the proof of Lemma \ref{Th:IndiscerniblesCofinalInSkolemHull}. 
\Square
\end{proof}

\begin{Lemma1} \label{Th:SkolemFunctionsFixedBelow2}
Let $\mathcal{M} \models \mathbf{W}_3$. For all $\mathcal{L}^*_{\hat{f}, \bar{\beta}}$ Skolem functions $\mathbf{t}(x_1, \ldots, x_n)$ and for all $i_1 < \ldots < i_n$ and $j_1 < \ldots < j_n$ in $\mathbb{Z}$, if for all $k \in \mathbb{Z}$, $\mathcal{M} \models \mathbf{t}(c_{i_1}, \ldots, c_{i_n}) < c_k$ then 
$$\mathcal{M} \models \mathbf{t}(c_{i_1}, \ldots, c_{i_n})= \mathbf{t}(c_{j_1}, \ldots, c_{j_n}).$$
\end{Lemma1}

\begin{proof}
Using the same arguments as the proof of Lemma \ref{Th:SkolemFunctionsFixedBelow}.
\Square
\end{proof}

\noindent This allows us to show that $\mathcal{L}_{\bar{\beta}}$-structures which are $\omega$-models of $\mathrm{Mac}$ and non-standard below $\bar{\beta}$ yield elementarily equivalent models admitting an automorphism that does not move any natural number down but moves an ordinal less than $\bar{\beta}$ down.

\begin{Theorems1} \label{Th:AutomorphismsFromNonStandardOmegaModels}
If $\mathcal{M}= \langle M^\mathcal{M}, \in^\mathcal{M}, \bar{\beta}^\mathcal{M} \rangle$ is an $\omega$-model of $\mathrm{Mac}$ such that\\ 
$\mathcal{M} \models \bar{\beta} \textrm{ is an ordinal}$ and $\mathcal{M}$ is non-standard below $\bar{\beta}^\mathcal{M}$ then there exists an $\mathcal{L}_{\bar{\beta}}$-structure $\mathcal{N} \equiv \mathcal{M}$ admitting an automorphism $j: \mathcal{N} \longrightarrow \mathcal{N}$ such that
\begin{itemize}
\item[(i)] $\mathcal{N} \models j(n) \geq n$ for all $n \in (\omega^\mathcal{N})^*$, 
\item[(ii)] there exists $\alpha \in (\bar{\beta}^\mathcal{N})^*$ such that $\mathcal{N} \models j(\alpha) < \alpha$. 
\end{itemize} 
\end{Theorems1}

\begin{proof}
Let $\mathcal{M}= \langle M^\mathcal{M}, \in^\mathcal{M}, \bar{\beta}^\mathcal{M} \rangle$ be an $\omega$-model of $\mathrm{Mac}$ such that 
$$\mathcal{M} \models \bar{\beta} \textrm{ is an ordinal}$$ 
and $\mathcal{M}$ is non-standard below $\bar{\beta}^\mathcal{M}$. By Lemma \ref{Th:ConsistencyOfW3} there exists an $\mathcal{L}^{\mathbf{ind}}_{\hat{f}, \bar{\beta}}$-structure $\mathcal{Q} \models \mathbf{W}_3$ such that $\mathcal{Q} \mid_{\mathcal{L}_{\bar{\beta}}} \equiv \mathcal{M}$. Let $C= \{ c_i^\mathcal{Q} \mid i \in \mathbb{Z} \}$. Let $\mathcal{N}= \mathcal{H}_{\mathcal{L}_{\hat{f}, \bar{\beta}}}^{\mathcal{Q}}(C)$. Therefore $\mathcal{N} \mid_{\mathcal{L}_{\bar{\beta}}} \equiv \mathcal{M}$. Define $j: \mathcal{N} \longrightarrow \mathcal{N}$ by 
$$j(c_i^\mathcal{Q})= c_{i+1}^\mathcal{Q} \textrm{ for all } i \in \mathbb{Z},$$
$$\textrm{for all } \mathcal{L}^*_{\hat{f}, \bar{\beta}} \textrm{ Skolem functions } \mathbf{t}(x_1, \ldots, x_n) \textrm{ and for all } i_1 < \ldots < i_n \textrm{ in } \mathbb{Z},$$
$$j(\mathbf{t}^\mathcal{Q}(c_{i_1}^\mathcal{Q}, \ldots, c_{i_n}^\mathcal{Q}))= \mathbf{t}^\mathcal{Q}(j(c_{i_1}^\mathcal{Q}), \ldots, j(c_{i_n}^\mathcal{Q})).$$
Now, $j$ is an automorphism of the $\mathcal{L}_{\hat{f}, \bar{\beta}}$-structure $\mathcal{N}$. The same arguments used in the proof of Theorem \ref{Th:AutomorphismMovingNoNaturalDown} reveal that if $\mathbf{t}(x_1, \ldots, x_n)$ is an $\mathcal{L}^*_{\hat{f}, \bar{\beta}}$ Skolem function and $i_1 < \ldots < i_n$ are in $\mathbb{Z}$ such that $\mathcal{N} \models \mathbf{t}(c_{i_1}^\mathcal{Q}, \ldots, c_{i_n}^\mathcal{Q}) \in \omega$ then 
$$\mathcal{N} \models j(\mathbf{t}(c_{i_1}^\mathcal{Q}, \ldots, c_{i_n}^\mathcal{Q})) \geq \mathbf{t}(c_{i_1}^\mathcal{Q}, \ldots, c_{i_n}^\mathcal{Q}).$$
Now, $j(\hat{f}^{\mathcal{Q}}(c_0^\mathcal{Q}))= \hat{f}^\mathcal{Q}(c_1^\mathcal{Q}) \in (\bar{\beta}^\mathcal{Q})^*$, and
$$\mathcal{Q} \models \hat{f}(c_1) < \hat{f}(c_0) < \bar{\beta}.$$
Therefore, 
$$\mathcal{N} \models j(\hat{f}(c_0^\mathcal{Q})) < \hat{f}(c_0^\mathcal{Q}) \textrm{ and } \hat{f}^\mathcal{N}(c_0^\mathcal{Q}) \in (\bar{\beta}^\mathcal{N})^*.$$
Therefore $j$ is the desired automorphism of the model $\mathcal{N} \mid_{\mathcal{L}_{\bar{\beta}}}$. 
\Square
\end{proof} 

Theorem \ref{Th:AutomorphismsFromNonStandardOmegaModels} allows us to give a positive answer to Question \ref{Q:HolmesQuestion} by showing that every complete consistent extension of $\mathrm{ZFC}$ has a model admitting an automorphism that does not move any natural number down but does move an ordinal less than the first non-recursive ordinal down. Using the same argument that we used in the proof of Theorem \ref{Th:ZFCwithAutThatDoesNotMoveAnyOrdinalDown} we will consider a model $\mathcal{M}$ of a complete consistent extension of $\mathrm{ZFC}$ with a limit ordinal $\lambda \in M^\mathcal{M}$ such that externally $\langle V_\lambda^\mathcal{M}, \in^\mathcal{M} \rangle$ is elementarily equivalent with $\mathcal{M}$. Using an argument developed by Harvey Friedman in \cite{fri73} we will apply the Barwise Compactness Theorem (Theorem \ref{Th:BarwiseCompactness}) to $V_\lambda$ inside $\mathcal{M}$ to obtain an $\omega$-model of an extension of $\mathrm{Mac}$ that is non-standard below $\omega_1^{\mathrm{ck}}$. Applying Theorem \ref{Th:AutomorphismsFromNonStandardOmegaModels} to this non-standard $\omega$-model will yield the model with the desired automorphism.

\begin{Theorems1} \label{Th:HolmesAutomorphism}
If $T \supseteq \mathrm{ZFC}$ is a complete and consistent $\mathcal{L}$-theory then there is an $\mathcal{L}$-structure $\mathcal{N} \models T$ admitting an automorphism $j: \mathcal{N} \longrightarrow \mathcal{N}$ such that
\begin{itemize}
\item[(i)] $\mathcal{N} \models j(n) \geq n$ for all $n \in (\omega^\mathcal{N})^*$,
\item[(ii)] there exists an $\alpha \in ((\omega_1^{\mathrm{ck}})^\mathcal{N})^*$ such that $\mathcal{N} \models j(\alpha) < \alpha$.
\end{itemize}
\end{Theorems1}

\begin{proof}
Let $T \supseteq \mathrm{ZFC}$ be a complete and consistent $\mathcal{L}$-theory. Let $\bar{\lambda}$ be a new constant symbol. Define $T^\star \supseteq T$ to be the $\mathcal{L}_{\bar{\lambda}}$-theory that includes the axioms:
\begin{itemize}
\item[(I)] $\bar{\lambda} > \omega_1^{\mathrm{ck}}$ is a limit ordinal, 
\item[(II)] for all $\mathcal{L}$-sentences, $\sigma$,
$$(\langle V_{\bar{\lambda}}, \in \rangle \models \sigma) \iff \sigma.$$ 
\end{itemize} 
Let $\Delta \subseteq T^\star$ be finite. Let $\mathcal{M}^\prime$ be an $\mathcal{L}$-structure such that $\mathcal{M}^\prime \models T$. Using the Reflection Principle we can expand $\mathcal{M}^\prime$ to an $\mathcal{L}_{\bar{\lambda}}$-structure $\mathcal{M}^{\prime\prime}$ such that 
$$\mathcal{M}^{\prime\prime} \models \Delta.$$
Using the Compactness Theorem we can find an $\mathcal{L}_{\bar{\lambda}}$-structure $\mathcal{M} \models T^\star$ such that $\omega^\mathcal{M}$ is non-standard. Let $a \in (\omega^\mathcal{M})^*$ be non-standard.\\ 
We work inside $\mathcal{M}$. Let $A= L_{\omega_1^{\mathrm{ck}}}$. Let $\mathcal{L}^\prime$ be the extension of $\mathcal{L}$ obtained by adding
\begin{itemize}
\item[(i)] the constant symbol $\bar{\beta}$, 
\item[(ii)] constant symbols $\hat{\alpha}$ for every $\alpha \in \omega_1^{\mathrm{ck}}$,
\item[(iii)] constant symbol $\mathbf{c}$. 
\end{itemize}
Let $\Gamma= \{ \phi \in \mathbf{Th}_{\mathcal{L}_{\bar{\beta}}}(\langle V_{\bar{\lambda}}, \in, \omega_1^{\mathrm{ck}} \rangle) \mid \ulcorner \phi \urcorner \in V_a \}$.\\
Define the $(\mathcal{L}^\prime_{\omega_1 \omega})_A$-theory $Q$ with axioms
\begin{itemize}
\item[(i)] $\Gamma$, 
\item[(ii)] $\hat{\zeta} \in \hat{\nu}$ for all $\zeta \in \nu \in \omega_1^{\mathrm{ck}}$, 
\item[(iii)] for all $\nu \in \omega_1^{\mathrm{ck}}$,
$$\forall x \left( x \in \hat{\nu} \Rightarrow \bigvee_{\zeta \in \nu} (x= \hat{\zeta})\right),$$
\item[(iv)] $\mathbf{c}$ is an ordinal,
\item[(v)] $\hat{\nu} \in \mathbf{c}$ for all $\nu \in \omega_1^{\mathrm{ck}}$,
\item[(vi)] $\forall x (x \in \mathbf{c} \cup \{\mathbf{c}\} \Rightarrow (x \textrm{ is not admissible}))$.
\end{itemize}
Since $\Gamma$ is finite, $Q$ is $\Sigma_1(\mathcal{L})$ over $A$. Let $Q^\prime \subseteq Q$ be such that $Q^\prime \in A$. Therefore
$$\{ \alpha \in \omega_1^{\mathrm{ck}} \mid \hat{\alpha} \textrm{ is mentioned in } Q^\prime \} \textrm{ is bounded in } \omega_1^{\mathrm{ck}}.$$
Therefore we can find a $\xi \in \omega_1^{\mathrm{ck}}$ which is greater than every $\alpha \in \omega_1^{\mathrm{ck}}$ such that $\hat{\alpha}$ is mentioned in $Q^\prime$. We expand $\langle V_{\bar{\lambda}}, \in, \omega_1^{\mathrm{ck}} \rangle$ by
\begin{itemize}
\item[(i)] interpreting $\mathbf{c}$ using $\xi$, 
\item[(ii)] if $\gamma \leq \alpha$ and $\hat{\alpha}$ is mentioned in $Q^\prime$ then we interpret $\hat{\gamma}$ using $\gamma$, 
\end{itemize}
to obtain a structure $\mathcal{Q} \models Q^\prime$. Therefore, by the Barwise Compactness Theorem, the $(\mathcal{L}^\prime_{\omega_1 \omega})_A$-theory $Q$ is consistent. Let $\mathcal{N}^{\prime\prime} \models Q$ and let 
$\mathcal{N}^\prime= \langle (N^\prime)^{\mathcal{N}^\prime}, \in^{\mathcal{N}^\prime}, \bar{\beta}^{\mathcal{N}^\prime} \rangle$ be the $\mathcal{L}_{\bar{\beta}}$-reduct of $\mathcal{N}^{\prime\prime}$. The finite axiom scheme (i) ensures that $\mathcal{N}^\prime \models \mathrm{Mac}$. Axiom schemes (ii) and (iii) of $Q$ ensure that $\mathcal{N}^\prime$ is an $\omega$-model. Now, suppose that $\mathcal{N}^\prime$ is well-founded below $\bar{\beta}^{\mathcal{N}^\prime}$. Let $\eta$ be the least element of 
$$\{ \xi \in (N^{\prime})^{\mathcal{N}^\prime} \mid (\forall \nu \in \omega_1^{\mathrm{ck}}) \mathcal{N}^{\prime\prime} \models \hat{\nu} \in \xi \} .$$
Now, $\eta= \omega_1^{\mathrm{ck}}$ and $L_{\omega_1^{\mathrm{ck}}}= L_\eta^{\mathcal{N}^\prime}$. Therefore
$$\mathcal{N}^\prime \models \eta \textrm{ is admissible}.$$
But this contradicts axiom (vi) of $Q$ since
$$\mathcal{N}^{\prime\prime} \models \eta \leq \mathbf{c}.$$
Therefore $\mathcal{N}^\prime$ is an $\omega$-model of $\mathrm{Mac}$ that is non-standard below $\bar{\beta}^{\mathcal{N}^\prime}$. And by (i) of $Q$,
$$\mathcal{N}^\prime \models \Gamma.$$
Therefore by Theorem \ref{Th:AutomorphismsFromNonStandardOmegaModels} there exists an $\mathcal{L}_{\bar{\beta}}$-structure $\mathcal{N} \equiv \mathcal{N}^\prime$ admitting an automorphism $j: \mathcal{N} \longrightarrow \mathcal{N}$ such that
\begin{itemize}
\item[(i)] $\mathcal{N} \models j(n) \geq n$ for all $n \in (\omega^\mathcal{N})^*$, 
\item[(ii)] there exists $\alpha \in (\bar{\beta}^\mathcal{N})^*$ such that $\mathcal{N} \models j(\alpha) < \alpha$. 
\end{itemize}
Now, $\mathcal{N} \models \bar{\beta}= \omega_1^{\mathrm{ck}}$ and $\mathcal{N} \models \Gamma$. Therefore the ambient arithmetic believes that $\mathcal{N} \models T$. And $j: \mathcal{N} \mid_{\mathcal{L}} \longrightarrow \mathcal{N} \mid_{\mathcal{L}}$ is an automorphism such that
\begin{itemize}
\item[(i)] $\mathcal{N} \mid_{\mathcal{L}} \models j(n) \geq n$ for all $n \in (\omega^{\mathcal{N} \mid_{\mathcal{L}}})^*$,
\item[(ii)] there exists an $\alpha \in ((\omega_1^{\mathrm{ck}})^{\mathcal{N} \mid_{\mathcal{L}}})^*$ such that $\mathcal{N} \mid_{\mathcal{L}} \models j(\alpha) < \alpha$.
\end{itemize}            
\Square
\end{proof}

\section[NFU]{NFU} \label{Sec:NFU}

In this section we give a brief description of the modification of Quine's `New Foundations' set theory $\mathrm{NFU}$ that was first introduced by Ronald Jensen in \cite{jen69}. A detailed introduction to $\mathrm{NFU}$ which includes the specifics of how mathematics can be formalised in this theory can be found in the introductory textbook \cite{hol98}. We would also like to direct readers to \cite{solXX} which contains a more condensed introduction to $\mathrm{NFU}$.\\
\\
\indent By weakening the extensionality axiom in Quine's $\mathrm{NF}$ to allow objects that are not sets into the domain of discourse, Jensen \cite{jen69} defines the subsystem $\mathrm{NFU}$ of $\mathrm{NF}$ which he is able to show is equiconsistent with a weak subsystem of $\mathrm{ZFC}$. In contrast to $\mathrm{NF}$ (see \cite{spe53}), Jensen's modified system is consistent with both the Axiom of Choice and the negation of the Axiom of Infinity. As we have already mentioned, in this paper we will use $\mathrm{NFU}$ to denote Jensen's weakening of Quine's $\mathrm{NF}$ fortified with both the Axiom of Infinity and the Axiom of Choice. Following \cite{hol98} and \cite{solXX} the Axiom of Infinity will be obtained in $\mathrm{NFU}$ by introducing a type-level ordered pair as a primitive notion. The underlying language of $\mathrm{NFU}$ is the language of set theory ($\mathcal{L}$) endowed with a unary predicate $\mathcal{S}$ and a ternary predicate $P$. The predicate $\mathcal{S}$ will be used to distinguish sets from non-sets. The ternary predicate $P$ will define a pairing function. Before presenting the axioms of $\mathrm{NFU}$ we first need to define the notion of a stratified formula.

\begin{Definitions1} \label{Df:StratifiedFormulae}
An $\mathcal{L}_{\mathcal{S}, P}$-formula $\phi$ is stratified if and only if there is a function $\sigma$ from the variables appearing in $\phi$ to $\mathbb{N}$ such that 
\begin{itemize}
\item[(i)] if $\textrm{`}x=y\textrm{'}$ is a sub-formula of $\phi$ then $\sigma(\textrm{`}x\textrm{'})= \sigma(\textrm{`}y\textrm{'})$, 
\item[(ii)] if $\textrm{`}x \in y\textrm{'}$ is a sub-formula of $\phi$ then $\sigma(\textrm{`}y\textrm{'})= \sigma(\textrm{`}x\textrm{'}) + 1$, 
\item[(iii)] if $\textrm{`}P(x, y, z)\textrm{'}$ is a sub-formula of $\phi$ then $\sigma(\textrm{`}z\textrm{'})= \sigma(\textrm{`}x\textrm{'})= \sigma(\textrm{`}y\textrm{'})$. 
\end{itemize} 
\end{Definitions1}

\noindent Equipped with this definition we are able to give the axioms of $\mathrm{NFU}$.

\begin{Definitions1} \label{Df:NFU}
We define $\mathrm{NFU}$ to be the $\mathcal{L}_{\mathcal{S}, P}$-theory with axioms
\begin{itemize}
\item[](Extensionality for Sets) 
$$\forall x \forall y(\mathcal{S}(x) \land \mathcal{S}(y) \Rightarrow (x=y \iff \forall z(z \in x \iff z \in y))),$$ 
\item[](Stratified Comprehension) for all stratified $\mathcal{L}_{\mathcal{S}, P}$-formulae $\phi(x, \vec{z})$,
$$\forall \vec{z} \exists y (\mathcal{S}(y) \land \forall x(x \in y \iff \phi(x, \vec{z}))),$$
\item[](Axioms of Pairing) 
$$\forall x \forall y \exists z (P(x, y, z) \land \forall w (P(x, y, w) \Rightarrow w= z)),$$
$$\forall x \forall y \forall u \forall v \forall z (P(x, y, z) \land P(u, v, z) \Rightarrow u = x \land v = y),$$
\item[](Axiom of Choice) every set is well ordered.
\end{itemize}
\end{Definitions1}    

Stratified comprehension proves the existence of a universal set which we will denote $V$. As is suggested by the above axiomatisation, ordered pairs in $\mathrm{NFU}$ are coded using the predicate $P$. That is to say, we define
$$\langle x, y \rangle= z \textrm{ if and only if } P(x, y, z).$$
The map $x \mapsto \langle x, x \rangle$ is a set by Stratified Comprehension, and the Axiom of Pairing implies that this map is an injection that is not a surjection. This shows that, as we mentioned above, the Axioms of Pairing imply the Axiom of Infinity.\\
\\
\indent The existence of `large' sets allows $\mathrm{NFU}$ to represent cardinal numbers as equivalence classes of equipollent sets. If $X$ is a set then we write $|X|$ for the cardinal to which $X$ belongs. In a similar fashion ordinals can be represented by equivalence classes of isomorphic well-orderings. If $R$ is a well-order then we will write $[R]$ for ordinal to which $R$ belongs. These representations allow the properties of being a cardinal number and being an ordinal number to be expressed using stratified formulae. This means that stratified comprehension proves the existence of a set of all cardinals, which we will denote $\mathrm{CN}$, and the existence of a set of all ordinals, which we will denote $\mathrm{ON}$. If $\alpha$ and $\beta$ are ordinals we say that $\alpha < \beta$ if and only if for every $R \in \alpha$ and $S \in \beta$, there exists an order isomorphism from $R$ into a proper initial segment of $S$. The order $<$ well-orders $\mathrm{ON}$ and we will use $\Omega$ to denote the ordinal to which $<$ belongs. The fact that cardinals and ordinals are coded using disjoint classes of objects endows importance to the following definition which mirrors the definition of initial ordinal in $\mathrm{ZFC}$. If $R$ is a binary relation then we use $\mathrm{Dom}(R)$ to denote the set $\{ x \mid \exists y (\langle x, y \rangle \in R \lor \langle y, x \rangle \in R)\}$.

\begin{Definitions1}
If $\kappa$ is a cardinal then we define $\mathrm{init}(\kappa)$ to be the least ordinal $\alpha$ such that there exists an $R \in \alpha$ and an $X \in \kappa$ such that $\mathrm{Dom}(R)= X$.
\end{Definitions1} 

We will use $\mathbb{N}$ to denote the set of all finite cardinal numbers and $\mathrm{CNI}$ to denote the set of all infinite cardinal numbers. If $\kappa$ and $\lambda$ are cardinal numbers then we say that $\kappa < \lambda$ if and only if for all $X \in \kappa$ and $Y \in \lambda$ there is an injection from $X$ into $Y$ but no bijection between $X$ and $Y$. The axiom of choice ensures that this ordering is a well-ordering. We will use $\aleph_0$ to denote the least infinite cardinal number and $\aleph_1$ to denote the least uncountable cardinal. If $X \subseteq \mathrm{CN}$ then we use $\sup X$ to denote the $<$-least $\kappa$ such that for all $\lambda \in X$, $\lambda \leq \kappa$. It should be noted that the restriction of comprehension to stratified formulae prevents $\mathrm{NFU}$ from proving the existence of von Neumann ordinals corresponding to well-orderings (see \cite{hen69} and \cite{fh09}).\\ 
\\
\indent Another unorthodox feature of $\mathrm{NFU}$ is the fact that it fails to prove that an arbitrary set is the same size as its own set of singletons. We will use $\iota$ to denote the singleton operation.  

\begin{Definitions1}
We say that a set $x$ is Cantorian if and only if $x$ is the same size as $\iota``x$. We say that $x$ is strongly Cantorian if and only if $\iota \upharpoonright x$ is a set.  
\end{Definitions1}    

\noindent For example, if $x$ is a set that has size a concrete natural number then $\mathrm{NFU}$ proves that $x$ is strongly Cantorian. On the other hand, $\mathrm{NFU}$ proves that $V$ is not Cantorian. This behaviour motivates the definition of the $T$ operation on cardinals, ordinals and equivalence classes of well-founded relations.

\begin{Definitions1}
We define the $T$ operation on cardinals and equivalence classes of well-founded relations as follows:
\begin{itemize}
\item[(i)] if $\kappa$ is a cardinal then
$$T(\kappa)= |\iota``A| \textrm{ where } A \in \kappa,$$  
\item[(ii)] if $[W]$ is the equivalence class of isomorphic well-founded relations with $W \in [W]$ then
$$T([W])= \textrm{ the class of relations isomorphic to } \{ \langle \iota x, \iota y \rangle \mid \langle x, y \rangle \in W\}.$$ 
\end{itemize}
\end{Definitions1}

\noindent It should be noted that any stratification of the formula $x= T(y)$ assigns a type to $x$ that is one higher than the type of $y$. In fact, if $\mathrm{NFU}$ is consistent then the collection $\{ \langle x, y \rangle \mid x= T(y) \}$ is a proper class. Another definable set operation that is inhomogeneous in the sense that any stratification of the formula defining this operation assigns different types to the arguments and the result is the powerset operation ($\mathcal{P}$). This inhomogeneity prevents $\mathrm{NFU}$ from proving Cantor's Theorem, and indeed the cardinal $|V|$ is an explicit counterexample to this result since $\mathcal{P}(V) \subseteq V$. Despite this failure, $\mathrm{NFU}$ is able to prove that for all sets $X$, $T(|X|) < |\mathcal{P}(X)|$. The following definition facilitates a type-homogeneous notion of cardinal exponentiation.

\begin{Definitions1}
If $\kappa$ is a cardinal then we define
$$2^{\kappa}= \left\{\begin{array}{ll}
\lambda & \textrm{if there exists an } X \in \kappa \textrm{ and a cardinal } \lambda \textrm{ s.t. } T(\lambda)= |\mathcal{P}(X)| \\
\emptyset & \textrm{otherwise}
\end{array} \right)$$  
\end{Definitions1}

\noindent The operation that sends $\kappa \mapsto 2^{\kappa}$ is definable by an $\mathcal{L}_{\mathcal{S}, P}$-formula that admits a stratification which assigns the same type to $\kappa$ and $2^{\kappa}$. The following result shows that this exponentiation operation possesses the strictly inflationary property that we intuitively associate with cardinal exponentiation. 

\begin{Lemma1} \label{Th:ExpInflationary}
Let $\kappa$ be a cardinal. If $2^{\kappa} \neq \emptyset$ then $\kappa < 2^{\kappa}$. \Square
\end{Lemma1}

\noindent This allows us to define an analogue of the Beth operation.

\begin{Definitions1}
We define {\small
$$\bar{\bar{\beth}}= \bigcap \left\{ X \subseteq \mathrm{CNI} \mid (\aleph_0 \in X) \land (\forall \kappa \in \mathrm{CNI})(\kappa \in X \Rightarrow 2^\kappa \in X) \land \forall Y ( Y \subseteq X \Rightarrow \sup Y \in X)\right\}.$$}
\end{Definitions1}

\noindent Stratified Comprehension ensures that $\bar{\bar{\beth}}$ is a set and Lemma \ref{Th:ExpInflationary} implies that $\bar{\bar{\beth}}$ is well-ordered by the natural ordering on cardinal numbers.

\begin{Definitions1}
If $\alpha$ is ordinal such that $\alpha < [< \cap(\bar{\bar{\beth}} \times \bar{\bar{\beth}})]$ then we define $\beth_\alpha^{TT}$ to be the $\alpha^{\textrm{th}}$ member of $\bar{\bar{\beth}}$.
\end{Definitions1}

The $T$ operation allows us to make explicit the relationship between an ordinal $\alpha$ and the initial segment of $\mathrm{ON}$ defined by $\alpha$.

\begin{Lemma1} \label{Th:OrderTypeOfOrdinalsLessThan}
If $\alpha$ is an ordinal and $X= \{ \beta \in \mathrm{ON} \mid \beta < \alpha \}$ then $[< \cap X \times X]= T^2(\alpha)$. \Square  
\end{Lemma1} 

Asserting that the $T$ operation is well-behaved is a natural way of extending the axioms of $\mathrm{NFU}$. One family of examples of this type of extension is Rosser's Axiom of Counting \cite{ros53} and two weakenings introduced in \cite{for77}.
\begin{itemize}
\item[]($\mathrm{AxCount}$) $(\forall n \in \mathbb{N}) T(n)=n$, 
\item[]($\mathrm{AxCount}_\leq$) $(\forall n \in \mathbb{N}) n \leq T(n)$,
\item[]($\mathrm{AxCount}_\geq$) $(\forall n \in \mathbb{N}) n \geq T(n)$. 
\end{itemize}

\noindent In \cite{ros53} it is shown that $\mathrm{AxCount}$ facilitates induction for unstratified properties. This added strength allows $\mathrm{NFU}+\mathrm{AxCount}$ to prove the consistency of Zermelo set theory \cite{hin75} and implies that $\mathrm{NFU}+\mathrm{AxCount}$ is consistent with the existence of classical sets such as the von Neumann $\omega$ and $V_\omega$ \cite{for06}. In \cite{for77}, \cite{hin75} and \cite{for06} it is observed that some of these consequences of $\mathrm{AxCount}$ also follow from the weakening $\mathrm{AxCount}_\leq$.\\   
\\
\indent In \cite{hol98} Holmes adapts the techniques developed by Rolland Hinnion in \cite{hin75} to show that fragments of Zermelo-Fraenkel set theory can be interpreted in classes of topped well-founded extensional relations in $\mathrm{NFU}$. 

\begin{Definitions1}
Let $R$ be a well-founded relation. For all $x \in \mathrm{Dom}(R)$, define $\mathrm{seg}_R(x)$ by induction
$$S_0= \{ \langle z, x \rangle \mid \langle z, x \rangle \in R \},$$
$$S_{n+1}= \{ \langle w, z \rangle \mid z \in \mathrm{Dom}(S_n) \land \langle w, z \rangle \in R \},$$
$$\mathrm{seg}_R(x)= \bigcup_{n \in \mathbb{N}} S_n.$$
\end{Definitions1}

\begin{Definitions1} \label{Df:BFEXT}
We say that a relation $R$ is a topped well-founded extensional relation (BFEXT) if and only if $R$ is well-founded extensional and there exists an $x \in \mathrm{Dom}(R)$ such that $R= \mathrm{seg}_R(x)$.  
\end{Definitions1}

\noindent Well-foundedness implies that if $R$ is a BFEXT then the $x \in \mathrm{Dom}(R)$ such that $R= \mathrm{seg}_R(x)$ is unique--- we will use $\mathbf{t}_R$ to denote this element.

\begin{Definitions1} \label{Df:EqualityAndMembershipForBFEXTs}
Let $R$ and $S$ be BFEXTs. We say that $R \cong S$ if and only if there is a bijection $f: \mathrm{Dom}(R) \longrightarrow \mathrm{Dom}(S)$ such that 
$$\langle x, y \rangle \in R \textrm{ if and only if } \langle f(x), f(y) \rangle \in S.$$
We say that $R \epsilon S$ if and only if there exists an $x \in \mathrm{Dom}(S)$ such that $\langle x, \mathbf{t}_S \rangle \in S$ and $R \cong \mathrm{seg}_S(x)$. 
\end{Definitions1}

If $R$ is a BFEXT then stratified comprehension allows $\mathrm{NFU}$ to prove that the following collections are sets:
\begin{equation} \label{Df:BFEXTEquivRelationInNFU}
[R]= \{ S \mid (S \textrm{ is a BFEXT }) \land (R \cong S) \}
\end{equation}
\begin{equation} \label{Df:BFinNFU}
\mathrm{BF}= \{ [R] \mid R \textrm{ is a BFEXT } \}
\end{equation}
We abuse notation and define the relation $\epsilon$ on $\mathrm{BF}$. For all $X, Y \in \mathrm{BF}$ we define
\begin{equation}
X \epsilon Y \textrm{ if and only if there exists } R \in X \textrm{ and } S \in Y \textrm{ such that } R \epsilon S.
\end{equation}

In \cite{solXX} Robert Solovay observes, without proof, that the structure $\langle \mathrm{BF}, \epsilon \rangle$ looks like $\langle H_{\kappa^+}, \in \rangle$ in the sense that there is a cardinal $\kappa$ such that
\begin{itemize}
\item[(I)] $\langle \mathrm{BF}, \epsilon \rangle$ is well-founded and extensional, 
\item[(II)] for every BFEXT $R$,
$$|\{ X \in \mathrm{BF} \mid X \epsilon [R] \}| \leq \kappa,$$ 
\item[(III)] if $X \subseteq \mathrm{BF}$ with $|X| \leq \kappa$ then there exists a BFEXT $R$ such that
$$X = \{ Y \mid Y \epsilon [R] \}.$$ 
\end{itemize}
Property (I) above is proved by Hinnion in \cite{hin75}:

\begin{Theorems1}
(Hinnion \cite{hin75}) The theory $\mathrm{NFU}$ proves that the structure $\langle \mathrm{BF}, \epsilon \rangle$ is well-founded and extensional. \Square
\end{Theorems1}

For the sake of completeness we will prove that if $\kappa= T^2(|V|)$ then properties (II) and (III) above hold. 

\begin{Theorems1} \label{Th:ExtensionsOfBFEXTsAreLessThanT2V}
The theory $\mathrm{NFU}$ proves that if $R$ is a BFEXT then
$$|\{ X \in \mathrm{BF} \mid X \epsilon [R] \}| \leq T^2(|V|).$$  
\end{Theorems1}

\begin{proof}
We work in $\mathrm{NFU}$. Let $R$ be a BFEXT. Let $\lhd$ be a well-ordering of $V$. Let $Y= \{ X \in \mathrm{BF} \mid X \epsilon [R] \}$. Define $f: Y \longrightarrow \iota^2``V$ by
\begin{equation} \label{Df:EmbeddingOfBFEXTensionIntoIotaV}
\begin{array}{lcl}
f(X)= \{\{x\}\} & \textrm{ if and only if } & S^{\prime\prime} \textrm{ is the } \lhd \textrm{-least element of } [R] \textrm{ and}\\
& & \textrm{there exists } S^\prime \in X \textrm{ and } f: S^\prime \longrightarrow S^{\prime\prime} \textrm{ witnessing}\\
& & \textrm{the fact that } S^\prime \cong \mathrm{seg}_{S^{\prime \prime}}(x) \textrm{ and } \langle x, \mathbf{t}_{S^{\prime\prime}} \rangle \in S^{\prime \prime}. 
\end{array}
\end{equation}
It is clear from (\ref{Df:EmbeddingOfBFEXTensionIntoIotaV}) that $f$ is defined by a stratified formula. Therefore $\mathrm{NFU}$ proves that $f$ is a set. The fact that $[R]$ is an equivalence class of isomorphic BFEXTs implies that $f$ is an injective function. Therefore $|Y| \leq T^2(|V|)$. 
\Square
\end{proof}

An important tool for building BFEXTs is the following Lemma, proved in \cite{hol98}, which shows that in $\mathrm{NFU}$ well-founded relations can be collapsed onto BFEXTs.

\begin{Lemma1} \label{Df:CollapsingLemmaForBFEXTs}
(Collapsing Lemma \cite{hol98}) The theory $\mathrm{NFU}$ proves that if $R$ is a well-founded relation then there exists an extensional well-founded relation $S$ and a map $g: \mathrm{Dom}(R) \longrightarrow \mathrm{Dom}(S)$ that is onto $\mathrm{Dom}(S)$ such that for all $x \in \mathrm{Dom}(R)$ and $z \in \mathrm{Dom}(S)$, $\langle z, g(x) \rangle \in S$ if and only if there exists a $y \in \mathrm{Dom}(R)$ such that $z= g(y)$ and $\langle y, x \rangle \in R$. Moreover, the relation $S$ is uniquely determined up to isomorphism by $R$. \Square   
\end{Lemma1}

This allows us to prove that if $\kappa= T^2(|V|)$ then property (III) above holds.

\begin{Theorems1} \label{Th:SetOfBFEXTsLessThanT2VAreExtensions}
The theory $\mathrm{NFU}$ proves that if $X \subseteq \mathrm{BF}$ is such that $|X| \leq T^2(|V|)$ then there exists a BFEXT $W$ such that 
$$X= \{ Y \in \mathrm{BF} \mid Y \epsilon [W] \}.$$ 
\end{Theorems1}

\begin{proof}
We work in $\mathrm{NFU}$. Let $\lhd$ be a well-ordering of $V$. Let $X \subseteq \mathrm{BF}$ be such that $|X| \leq T^2(|V|)$. Let $f: X \longrightarrow \iota^2``V$ be an injection. Let
$$D= \{ \mathrm{Dom}(R) \times \{a \} \mid (R \textrm{ is the } \lhd \textrm{-least element of }[R])\land ([R] \in X) \land(f([R])= \{\{a\}\})\}.$$
Let $\xi$ be such that $\xi \notin \bigcup D$. Define $S \subseteq \bigcup D \times \bigcup D$ by {\small
$$\begin{array}{lcl}
\langle x, y \rangle \in S & \textrm{ if and only if } & x= \langle w, a \rangle \textrm{ and } y= \langle z, a \rangle \textrm{ and } (R \textrm{ is the } \lhd \textrm{-least element of } [R])\\
& & \textrm{and }([R] \in X) \textrm{ and } (f([R])= \{\{a\}\}) \textrm{ and } (\langle w, z \rangle \in R) \textrm{ or},\\
& & x= \langle w, a \rangle \textrm{ and } y= \xi \textrm{ and } (R \textrm{ is the } \lhd \textrm{-least element of } [R])\\
& & \textrm{and } ([R] \in X) \textrm{ and } (f([R])= \{\{a\}\}) \textrm{ and } w= \mathbf{t}_S.
\end{array}$$}
Now, $S$ is a well-founded relation with a unique top $\xi$. Moreover, for all BFEXTs $R$, if $[R] \in X$ then there exists an $x \in \mathrm{Dom}(S)$ such that $\langle x, \xi \rangle \in S$ and $R \cong \mathrm{seg}_S(x)$. Let $W$ be the extensional collapse of $S$ obtained using Lemma \ref{Df:CollapsingLemmaForBFEXTs}. Now, $W$ is a BFEXT and 
$$X= \{ Y \in \mathrm{BF} \mid Y \epsilon [W] \}.$$   
\Square
\end{proof}

This characterisation implies that the structure $\langle \mathrm{BF}, \epsilon \rangle$ is a model of $\mathrm{ZFC}^-$.

\begin{Coroll1} \label{Th:BFEXTsModelZFCminus}
(Hinnion \cite{hin75}, Holmes \cite{hol98}) The theory $\mathrm{NFU}$ proves that
$$\langle \mathrm{BF}, \epsilon \rangle \models \mathrm{ZFC}^-.$$
\Square
\end{Coroll1}

The operation $T$ restricted to the set $\mathrm{BF}$ is an embedding of the structure $\langle \mathrm{BF}, \epsilon \rangle$ into itself. That is to say, for all $[R], [S] \in \mathrm{BF}$,
$$\textrm{both } T([R]) \textrm{ and } T([S]) \textrm{ are BFEXTs and}$$
$$[R] \epsilon [S] \textrm{ if and only if } T([R]) \epsilon T([S]).$$
Hinnion \cite{hin75} noted that the set $\mathrm{ON}$ can be bijectively mapped onto the von Neumann ordinals of the Zermelian structure $\langle \mathrm{BF}, \epsilon \rangle$. Define $\mathbf{k}: \mathrm{ON} \longrightarrow \mathrm{BF}$ by
\begin{equation} \label{Df:EmbeddingK}
\mathbf{k}([R])= [R \cup \{ \langle x, x^* \rangle \mid x \in \mathrm{Dom}(R) \}] \textrm{ where } R \textrm{ is well-order such that } x^* \notin \mathrm{Dom}(R).
\end{equation}
The definition of $\mathbf{k}$ admits a stratification that assigns the same type to both the result and the argument. Therefore stratified comprehension can be used to show that the map $\mathbf{k}: \mathrm{ON} \longrightarrow \mathrm{Ord}^{\langle \mathrm{BF}, \epsilon \rangle}$ is a set. It is also clear that the map $\mathbf{k}$ is order preserving and for all equivalence classes of well-orderings $[R]$,
$$\mathbf{k}(T([R]))= T(\mathbf{k}([R])).$$
We can extend the notions of Cantorian and strongly Cantorian to ordinals and equivalence classes of BFEXTs.

\begin{Definitions1}
We say that an ordinal $\alpha$ is Cantorian if and only if $T(\alpha)= \alpha$. We say that $\alpha$ is strongly Cantorian if and only if for all $\beta \leq \alpha$, $T(\beta)= \beta$.
\end{Definitions1}

\begin{Definitions1} \label{Df:CandSCBFEXTs}
Let $R$ be a BFEXT. We say that $[R]$ is Cantorian if and only if $T([R])= [R]$. We say that $[R]$ is strongly Cantorian if and only if for all BFEXTs $S$ with 
$$\langle \mathrm{BF}, \epsilon \rangle \models [S] \epsilon \mathrm{TC}(\{[R]\}),$$
$T([S])= [S]$.
\end{Definitions1}

In \cite{jen69} Jensen proves an analogue of Specker's Theorem \cite{spe53} showing that models of $\mathrm{NFU}$ can be produced from models of a modification of Hao Wang's theory of types (see \cite{wan52}), obtained by weakening the extensionality axiom scheme, that admit a type-shifting automorphism.

\begin{Definitions1}
The language $\mathcal{L}_{\mathrm{T}\mathbb{Z}\mathrm{TU}}$ is the $\mathbb{Z}$-sorted theory with variables\\ 
$x_i^1, y_i^1, z_i^1, x_i^2, y_i^2, z_i^2, \ldots$ for each $i \in \mathbb{Z}$ and
\begin{itemize}
\item[(i)] binary predicates $\in_i$ for each $i \in \mathbb{Z}$,
\item[(ii)] binary predicates $=_i$ for each $i \in \mathbb{Z}$, 
\item[(iii)] unary predicates $\mathcal{S}_i$ for each $i \in \mathbb{Z}$,
\item[(iv)] ternary predicates $P_i$ for each $i \in \mathbb{Z}$.
\end{itemize}
Well-formed formulae of $\mathcal{L}_{\mathrm{T}\mathbb{Z}\mathrm{TU}}$ are built up inductively from the atomic formulae\\
$x_i^n \in_i y_{i+1}^m$, $x_i^n =_i y_i^m$, $\mathcal{S}_i(x_i^n)$ and $P_i(x_i^m, y_i^k, z_i^n)$ for all $i \in \mathbb{Z}$ and for all $n, m, k \in \mathbb{N}$. 
\end{Definitions1}

An $\mathcal{L}_{\mathrm{T}\mathbb{Z}\mathrm{TU}}$-structure $\mathcal{M}$ consists of domains $M_i^\mathcal{M}$ and relations $\in_i^\mathcal{M}$, $\mathcal{S}_i^{\mathcal{M}}$ and $P_i^{\mathcal{M}}$ for each $i \in \mathbb{Z}$. For each $i \in \mathbb{Z}$, the predicates $=_i^\mathcal{M}$ will be interpreted as true identity in the domain $M_i^\mathcal{M}$. Quantifiers appearing in well-formed $\mathcal{L}_{\mathrm{T}\mathbb{Z}\mathrm{TU}}$-formulae range over objects of the same type as variable directly to the right of the quantifier. For example $\forall x_i$ reads `for all objects of type $i$'. We are now in a position to describe the theory of types, which we will call $\mathrm{T}\mathbb{Z}\mathrm{TU}$, related to $\mathrm{NFU}$ introduced in \cite{jen69}. We differ from \cite{jen69} by including both the axiom of infinity (at each type $i \in \mathbb{Z}$) and the axiom of choice (at each type $i \in \mathbb{Z}$) in our axiomatisation of $\mathrm{T}\mathbb{Z}\mathrm{TU}$. As with our axiomatisation of $\mathrm{NFU}$ we ensure that the axiom of infinity holds at each type $i \in \mathbb{Z}$ by specifying that for all $i \in \mathbb{Z}$, the predicate $P_i$ defines a type-level pairing function. 

\begin{Definitions1}
We define $\mathrm{T}\mathbb{Z}\mathrm{TU}$ to be the $\mathcal{L}_{\mathrm{T}\mathbb{Z}\mathrm{TU}}$-theory with axioms
\begin{itemize}
\item[](Extensionality for Sets) for all $i \in \mathbb{Z}$,
$$\forall x_{i+1} \forall y_{i+1} \left( \begin{array}{c}
\mathcal{S}_{i+1}(x_{i+1})\land \mathcal{S}_{i+1}(y_{i+1}) \Rightarrow\\
(x_{i+1} =_{i+1} y_{i+1} \iff \forall z_i(z_i \in_i x_{i+1} \iff z_i \in_i y_{i+1}))
\end{array}\right),$$  
\item[](Comprehension) for all $i \in \mathbb{Z}$ and for all $\mathcal{L}_{\mathrm{T}\mathbb{Z}\mathrm{TU}}$-formula $\phi(x_i, \vec{z})$,  
$$\forall \vec{z} \exists y_{i+1}( \mathcal{S}_{i+1}(y_{i+1}) \land \forall x_i(x_i \in_i y_{i+1} \iff \phi(x_i, \vec{z}))),$$
\item[](Pairing) for all $i \in \mathbb{Z}$,
$$\forall x_i \forall y_i \exists z_i (P_i(x_i, y_i, z_i) \land \forall w_i (P_i(x_i, y_i, w_i) \Rightarrow w_i =_i z_i)),$$
$$\forall x_i \forall y_i \forall u_i \forall v_i \forall z_i (P_i(x_i, y_i, z_i) \land P_i(u_i, v_i, z_i) \Rightarrow u_i =_i x_i \land v_i =_i y_i).$$
\item[](Axiom of Choice) for all $i \in \mathbb{Z}$, the set of all objects of type $i$ have a well ordering [at type $i+1$]. 
\end{itemize}
\end{Definitions1}

\noindent The theory $\mathrm{T}\mathbb{Z}\mathrm{TU}$ is Wang's theory \cite{wan52} fortified with the axiom of infinity and axiom of choice (at each type $i \in \mathbb{Z}$), and with extensionality weakened (at each type $i \in \mathbb{Z}$) to allow urelements into the domain of discourse.  

\begin{Definitions1}
If $\phi$ is an $\mathcal{L}_{\mathrm{T}\mathbb{Z}\mathrm{TU}}$-formula then we define $\phi^\dagger$ to be the $\mathcal{L}_{\mathcal{S}, P}$-formula obtained by deleting the subscripts from the predicate and function symbols appearing in $\phi$. 
\end{Definitions1} 

\noindent Note that if $\phi$ is a well-formed $\mathcal{L}_{\mathrm{T}\mathbb{Z}\mathrm{TU}}$-formula then $\phi^\dagger$ is stratified.

\begin{Definitions1}
Let $\mathcal{M}= \langle (M_i^{\mathcal{M}})_{i \in \mathbb{Z}}, (\in_i^{\mathcal{M}})_{i \in \mathbb{Z}}, (\mathcal{S}_i^{\mathcal{M}})_{i \in \mathbb{Z}}, (P_i^{\mathcal{M}})_{i \in \mathbb{Z}} \rangle$ be an $\mathcal{L}_{\mathrm{T}\mathbb{Z}\mathrm{TU}}$-structure. A type-shifting automorphism of $\mathcal{M}$ is a bijection 
$$\sigma: \bigcup_{i \in \mathbb{Z}} M_i^{\mathcal{M}} \longrightarrow \bigcup_{i \in \mathbb{Z}} M_i^{\mathcal{M}} \textrm{ such that for all } i \in \mathbb{Z},$$
\begin{itemize}
\item[(i)] $\sigma \upharpoonright M_i^{\mathcal{M}}$ is a bijection onto $M_{i+1}^{\mathcal{M}}$, 
\item[(ii)] for all $x \in M_i^\mathcal{M}$ and for all $y \in M_{i+1}^{\mathcal{M}}$,
$$x \in_i^\mathcal{M} y \textrm{ if and only if } \sigma(x) \in_{i+1}^\mathcal{M} \sigma(y),$$
\item[(iii)] for all $x \in M_i^\mathcal{M}$,
$$\mathcal{S}_i(x) \textrm{ if and only if } \mathcal{S}_{i+1}(\sigma(x)),$$
\item[(iv)] for all $x, y, z \in M_i^\mathcal{M}$,
$$P_i(x, y, z) \textrm{ if and only if } P_{i+1}(\sigma(x), \sigma(y), \sigma(z)).$$
\end{itemize}  
\end{Definitions1}

Equipped with this definition we are now in a position to recall Jensen's analogue of Specker's Theorem that links models of $\mathrm{NFU}$ to models of $\mathrm{T}\mathbb{Z}\mathrm{TU}$ admitting a type-shifting automorphism. Let $\mathcal{M}= \langle (M_i^\mathcal{M})_{i \in \mathbb{Z}}, (\in_i^\mathcal{M})_{i \in \mathbb{Z}}, (\mathcal{S}_i^\mathcal{M})_{i \in \mathbb{Z}}, (P_i^\mathcal{M})_{i \in \mathbb{Z}} \rangle$ be an $\mathcal{L}_{\mathrm{T}\mathbb{Z}\mathrm{TU}}$-structure such that $\mathcal{M}\models \mathrm{T}\mathbb{Z}\mathrm{TU}$ with a type-shifting automorphism $\sigma$. Define for all $x, y \in M_0^{\mathcal{M}}$,
\begin{equation} \label{Df:NFUMembership}
x \in_\sigma y \textrm{ if and only if } x \in_0^\mathcal{M} \sigma(y).
\end{equation}             
The membership relation defined on $M_0^{\mathcal{M}}$ by (\ref{Df:NFUMembership}) yields an $\mathcal{L}_{\mathcal{S}, P}$-structure\\ 
$\langle M_0^{\mathcal{M}}, \in_\sigma, \mathcal{S}_0^\mathcal{M}, P_0^\mathcal{M} \rangle$ that models $\mathrm{NFU}$. 

\begin{Theorems1} \label{Th:NFUfromTZTU}
(Jensen) Let $\mathcal{M}= \langle (M_i^\mathcal{M})_{i \in \mathbb{Z}}, (\in_i^\mathcal{M})_{i \in \mathbb{Z}}, (\mathcal{S}_i^\mathcal{M})_{i \in \mathbb{Z}}, (P_i^\mathcal{M})_{i \in \mathbb{Z}} \rangle$ be an $\mathcal{L}_{\mathrm{T}\mathbb{Z}\mathrm{TU}}$-structure. If $\mathcal{M}\models  \mathrm{T}\mathbb{Z}\mathrm{TU}$ and $\sigma$ is a type-shifting automorphism of $\mathcal{M}$ then \mbox{$\langle M_0^{\mathcal{M}}, \in_\sigma, \mathcal{S}_0^\mathcal{M}, P_0^{\mathcal{M}} \rangle$} is a model of $\mathrm{NFU}$. Moreover, for all $\mathcal{L}_{\mathrm{T}\mathbb{Z}\mathrm{TU}}$-formulae $\phi(x_{t_1}^1, \ldots, x_{t_n}^n)$ and for all $a_1, \ldots, a_n \in M_0^\mathcal{M}$,
$$\mathcal{M} \models \phi(\sigma^{t_1}(a_1), \ldots, \sigma^{t_n}(a_n)) \textrm{ if and only if } \langle M_0^{\mathcal{M}}, \in_\sigma, \mathcal{S}_0^\mathcal{M}, P_0^{\mathcal{M}} \rangle \models \phi^\dagger(a_1, \ldots, a_n).$$
\Square 
\end{Theorems1}

Maurice Boffa \cite{bof88} observed that models of $\mathrm{NFU}$ arise naturally from models of set theory that admit automorphism. Let $\mathcal{M}= \langle M^\mathcal{M}, \in^\mathcal{M} \rangle$ be an $\mathcal{L}$-structure such that $\mathcal{M} \models \mathrm{Mac}$. Let $j: \mathcal{M} \longrightarrow \mathcal{M}$ be an automorphism and let $c \in M^\mathcal{M}$ be such that
\begin{itemize}
\item[(i)] $\mathcal{M} \models c \textrm{ is infinite}$, 
\item[(ii)] $\mathcal{M} \models c \cup \mathcal{P}(c) \subseteq j(c)$. 
\end{itemize}
There is a natural model of $\mathrm{T}\mathbb{Z}\mathrm{TU}$ that arises from $\mathcal{M}$, $j$ and $c$. For all $i \in \mathbb{Z}$, define
\begin{equation} \label{Df:TZTUFromMacWithAut1}
M_i= j^i(c) \textrm{ and } E_i= \in^\mathcal{M} \cap (j^i(c) \times \mathcal{P}^\mathcal{M}(j^i(c))) \textrm{ and } S_i= \mathcal{P}^\mathcal{M}(j^{i-1}(c)). 
\end{equation}
For all $i \in \mathbb{Z}$, we can use the Axiom of Choice in $\mathrm{Mac}$ to define a ternary relation $\bar{P}_i$ on $M_i$ such that
\begin{equation} \label{Df:TZTUFromMacWithAut2}
\mathcal{M}_j^c= \langle (M_i)_{i \in \mathbb{Z}}, (E_i)_{i \in \mathbb{Z}}, (S_i)_{i \in \mathbb{Z}}, (\bar{P}_i)_{i \in \mathbb{Z}} \rangle \models \mathrm{T}\mathbb{Z}\mathrm{TU}. 
\end{equation}
Now, it is clear from this construction that $j$ is a type-shifting automorphism of the $\mathcal{L}_{\mathrm{T}\mathbb{Z}\mathrm{TU}}$-structure $\mathcal{M}_j^c$. This yields the following result:

\begin{Theorems1} \label{Th:NFUFromMacWithAut}
Let $\mathcal{M}$ be an $\mathcal{L}$-structure such that $\mathcal{M} \models \mathrm{Mac}$. Let $j: \mathcal{M} \longrightarrow \mathcal{M}$ be an automorphism and let $c \in M^\mathcal{M}$ be such that
\begin{itemize}
\item[(i)] $\mathcal{M} \models c \textrm{ is infinite}$, 
\item[(ii)] $\mathcal{M} \models c \cup \mathcal{P}(c) \subseteq j(c)$. 
\end{itemize}
If $\mathcal{M}_j^c$ is defined by (\ref{Df:TZTUFromMacWithAut1}) and (\ref{Df:TZTUFromMacWithAut2}) and $\in_j$ is defined by (\ref{Df:NFUMembership}) from the type-shifting automorphism $j$ acting on $\mathcal{M}_j^c$ then $\langle M_0, \in_j, S_0, \bar{P}_0 \rangle \models \mathrm{NFU}$. Moreover, for all $\mathcal{L}_{\mathrm{T}\mathbb{Z}\mathrm{TU}}$-formulae $\phi(x_{t_1}^1, \ldots, x_{t_n}^n)$ and for all $a_1, \ldots, a_n \in M_0$,
$$\mathcal{M}_j^c \models \phi(j^{t_1}(a_1), \ldots, j^{t_n}(a_n)) \textrm{ if and only if } \langle M_0, \in_j, S_0, \bar{P}_0 \rangle \models \phi^\dagger(a_1, \ldots, a_n).$$
\Square
\end{Theorems1}

\noindent The results of \cite{bof88} reveal that condition (ii) of Theorem \ref{Th:NFUFromMacWithAut} can be weakened to
$$\mathcal{M} \models |c \cup \mathcal{P}(c) | \leq |j(c)|.$$

In \cite{jen69}, Jensen identifies a weak subsystem of $\mathrm{ZFC}$ that is equiconsistent with $\mathrm{NFU}$. It should be noted that the arguments appearing in \cite{jen69} which claim to establish this equiconsistency are incomplete. However, using results proved by Thomas Forster and Richard Kaye that are reported in \cite{mat01}, Jensen's model constructions in \cite{jen69} yield the following equiconsistency result:    

\begin{Theorems1}
(Jensen \cite{jen69}) $\mathrm{NFU}$ is equiconsistent with $\mathrm{Mac}$. 
\Square
\end{Theorems1}

\section[The strength of $\mathrm{NFU}+\mathrm{AxCount}_\geq$]{The strength of $\mathrm{NFU}+\mathrm{AxCount}_\geq$} \label{Sec:NFUPlusAxCountGEQ}

In this section we will use Theorem \ref{Th:AutomorphismMovingNoNaturalDown} to prove two results concerning the strength of the theory $\mathrm{NFU}+\mathrm{AxCount}_\geq$. We will first show that $\mathrm{NFU}+\mathrm{AxCount}_\geq$ is unable to prove that $\mathrm{CNI}$ is infinite. We will then use the fact that the theory $\mathrm{NFU}+\mathrm{AxCount}_\leq$ is able to build $\omega$-models of Zermelo set theory in the structure $\langle \mathrm{BF}, \epsilon \rangle$ to show that $\mathrm{NFU}+ \mathrm{AxCount}_\leq$ proves the consistency of $\mathrm{NFU}+\mathrm{AxCount}_\geq$.  

\begin{Theorems1} \label{Th:ModelOfNFUAxCountGEQBethOmegaDoesNotExist}
There exists an $\mathcal{L}_{\mathcal{S}, P}$-structure $\mathcal{M}$ such that $\mathcal{M} \models \mathrm{NFU} + \mathrm{AxCount}_\geq$ and 
$$\mathcal{M} \models \mathrm{CNI} \textrm{ is finite}.$$
\end{Theorems1}

\begin{proof}
We work in the theory $\mathrm{ZFC}+\mathrm{GCH}$. The $\mathcal{L}$-structure $\langle V_{\omega+\omega}, \in \rangle$ is an $\omega$-model of Zermelo set theory plus $\mathrm{TCo}$. By Theorem \ref{Th:AutomorphismMovingNoNaturalDown} there is an $\mathcal{L}$-structure $\mathcal{N} \equiv \langle V_{\omega+\omega}, \in \rangle$ admitting an automorphism $j: \mathcal{N} \longrightarrow \mathcal{N}$ such that
\begin{itemize}
\item[(i)] $\mathcal{N} \models j(n) \geq n$ for all $n \in (\omega^\mathcal{N})^*$, 
\item[(ii)] there exists an $n \in (\omega^\mathcal{N})^*$ with $\mathcal{N} \models j(n) > n$. 
\end{itemize}
Since $\mathcal{N} \equiv \langle V_{\omega+\omega}, \in \rangle$, the rank function $\rho$ is well-defined and 
$$\mathcal{N} \models \forall \alpha ((\alpha \textrm{ is an ordinal}) \Rightarrow (V_\alpha \textrm{ exists})).$$
Let $n \in (\omega^\mathcal{N})^*$ be such that $\mathcal{N} \models j(n) > n$. Recalling (\ref{Df:TZTUFromMacWithAut1}) and (\ref{Df:TZTUFromMacWithAut2}) we can use $j$ and $V_{\omega+n}$ to define an $\mathcal{L}_{\mathrm{T}\mathbb{Z}\mathrm{TU}}$-structure $\mathcal{M}^{V_{\omega+n}^{\mathcal{N}}}_j \models \mathrm{T}\mathbb{Z}\mathrm{TU}$ such that $j$ is a type-shifting automorphism of $\mathcal{M}^{V_{\omega+n}^{\mathcal{N}}}_j$. Using (\ref{Df:NFUMembership}) and Theorem \ref{Th:NFUFromMacWithAut} we can define $\in_j$, $S$ and $\bar{P}$ on $V_{\omega+n}^\mathcal{N}$ so as 
$$\mathcal{M}= \langle V_{\omega+n}^\mathcal{N}, \in_j, S, \bar{P} \rangle \models \mathrm{NFU}.$$  
By the properties of $j$, for all $n \in (\omega^\mathcal{N})^*$, 
$$\mathcal{M}^{V_{\omega+n}^{\mathcal{N}}}_j \models |\iota``n| \leq j(n).$$
Therefore, by Theorem \ref{Th:NFUFromMacWithAut}, for all $n \in (\omega^\mathcal{N})^*$,
$$\mathcal{M} \models T(n) \leq n.$$
Therefore $\mathcal{M} \models \mathrm{AxCount}_\geq$. Now, let $\sigma_i$ be the $\mathcal{L}_{\mathrm{T}\mathbb{Z}\mathrm{TU}}$-sentence that says: `the set of infinite cardinals [at level $i$] is finite'. It follows from the Generalised Continuum Hypothesis that for all $i \in \mathbb{Z}$,
$$\mathcal{M}^{V_{\omega+n}^{\mathcal{N}}}_j \models \sigma_i.$$
Therefore $\mathcal{M} \models \mathrm{CNI} \textrm{ is finite}$.   
\Square
\end{proof}

\noindent By observing that $\bar{\bar{\beth}} \subseteq \mathrm{CNI}$ we get the following result:

\begin{Coroll1}
There exists an $\mathcal{L}_{\mathcal{S}, P}$-structure $\mathcal{M}$ such that $\mathcal{M} \models \mathrm{NFU} + \mathrm{AxCount}_\geq$ and 
$$\mathcal{M} \models \beth_\omega^{TT} \textrm{ does not exist}.$$
\Square
\end{Coroll1} 

We now turn to showing that the theory $\mathrm{NFU}+\mathrm{AxCount}_\leq$ is strictly stronger than $\mathrm{NFU}+\mathrm{AxCount}_\geq$. In \cite{hin75} Roland Hinnion shows that the ordinal strength endowed upon $\mathrm{NF}$ by $\mathrm{AxCount}_\leq$ means that there exists a point in the structure $\langle \mathrm{BF}, \epsilon \rangle$ that $\langle \mathrm{BF}, \epsilon \rangle$ believes is a transitive model of Zermelo set theory. Using the characterisation of the structure $\langle \mathrm{BF}, \epsilon \rangle$ proved in section \ref{Sec:NFU} we will reprove Hinnion's result in the context of $\mathrm{NFU}$. We make the following definitions in $\mathrm{NFU}$:

\begin{Definitions1}
Let $X \in \mathrm{BF}$. We define
$$\mathrm{Pow}(X)= \{ Y \in \mathrm{BF} \mid \langle \mathrm{BF}, \epsilon \rangle \models Y \subseteq [R] \}.$$
\end{Definitions1}

\begin{Definitions1}
Let $X \in \mathrm{BF}$. We define
$$X^{\mathrm{ext}}= \{ Y \in \mathrm{BF} \mid \langle \mathrm{BF}, \epsilon \rangle \models Y \in [R] \}.$$
\end{Definitions1}

Corollary \ref{Th:BFEXTsModelZFCminus} implies that $\mathrm{NFU}$ proves that
$$\langle \mathrm{BF}, \epsilon \rangle \models V_\omega \textrm{ exists}.$$ 
Adding $\mathrm{AxCount}_\leq$ to $\mathrm{NFU}$ allows us to show that the structure $\langle \mathrm{BF}, \epsilon \rangle$ believes that $V_{\omega+\omega}$ exists.

\begin{Lemma1} \label{Th:SizeOfPow}
The theory $\mathrm{NFU}$ proves that if $X \in \mathrm{BF}$ then 
$$|\mathrm{Pow}(X)|= 2^{|X^{\mathrm{ext}}|}.$$
\end{Lemma1}

\begin{proof}
We work in $\mathrm{NFU}$. Let $X \in \mathrm{BF}$. Define $f: \iota``\mathrm{Pow}(X) \longrightarrow \mathcal{P}(X^{\mathrm{ext}})$ by
$$\begin{array}{lcl}
f(\{Y\})= U & \textrm{ if and only if } & \forall V(V \in U \iff \langle \mathrm{BF}, \epsilon \rangle \models V \in Y).
\end{array}$$
Stratified comprehension ensures that $f$ is a set. It is clear that $f$ is a bijective function. Therefore $T(|\mathrm{Pow}(X)|)= |\mathcal{P}(X^{\mathrm{ext}})|$, which proves the lemma.
\Square
\end{proof}

\noindent This allows us to show that the theory $\mathrm{NFU}+\mathrm{AxCount}_\leq$ proves that for all internal natural numbers $n$, the structure $\langle \mathrm{BF}, \epsilon \rangle$ believes that $V_{\omega+n}$ exists.

\begin{Lemma1} \label{Th:BethNLessThanT2V}
The theory $\mathrm{NFU}+\mathrm{AxCount}_\leq$ proves that for all $n \in \mathbb{N}$, $\beth_n^{TT} \leq T^2(|V|)$.  
\end{Lemma1}

\begin{proof}
We work in the theory $\mathrm{NFU}+\mathrm{AxCount}_\leq$. Let $n$ be least such that $T^2(|V|) < \beth_n^{TT} \leq |V|$. Therefore $T^4(|V|) < T^2(\beth_n^{TT}) \leq T^2(|V|)$. And so, $\beth_{T^2(n)}^{TT}= T^2(\beth_n^{TT}) < \beth_n^{TT}$. But this contradicts the fact that $T^2(n) \geq n$.
\Square
\end{proof}

\begin{Theorems1} \label{Th:AxCountLEQImpliesVOmegaPlusNExists}
The theory $\mathrm{NFU}+\mathrm{AxCount}_\leq$ proves that for all $n \in \mathbb{N}$,
\begin{itemize}
\item[(i)] $\langle \mathrm{BF}, \epsilon \rangle \models V_{\omega+\mathbf{k}(n)} \textrm{ exists}$, 
\item[(ii)] $|(V_{\omega+\mathbf{k}(n)}^{\langle \mathrm{BF}, \epsilon \rangle})^{\mathrm{ext}}|= \beth_n^{TT}$. 
\end{itemize}
\end{Theorems1}

\begin{proof}
We work in the theory $\mathrm{NFU}+\mathrm{AxCount}_\leq$. We prove the theorem by induction. Corollary \ref{Th:BFEXTsModelZFCminus} implies that 
\begin{itemize}
\item[(i)] $\langle \mathrm{BF}, \epsilon \rangle \models V_\omega \textrm{ exists}$,
\item[(ii)] $|(V_{\omega}^{\langle \mathrm{BF}, \epsilon \rangle})^{\mathrm{ext}}|= \aleph_0= \beth_0^{TT}$.
\end{itemize}
Suppose that the theorem holds for $n \in \mathbb{N}$. Now, by Lemma \ref{Th:SizeOfPow},
$$|\mathrm{Pow}(V_{\omega+\mathbf{k}(n)}^{\langle \mathrm{BF}, \epsilon \rangle})| = 2^{|(V_{\omega+\mathbf{k}(n)}^{\langle \mathrm{BF}, \epsilon \rangle})^{\mathrm{ext}}|}= 2^{\beth_n^{TT}}= \beth_{n+1}^{TT}.$$
Therefore, by Lemma \ref{Th:BethNLessThanT2V}, $|\mathrm{Pow}(V_{\omega+\mathbf{k}(n)}^{\langle \mathrm{BF}, \epsilon \rangle})| \leq T^2(|V|)$. And so, by Theorem \ref{Th:SetOfBFEXTsLessThanT2VAreExtensions}, there exists an $X \in \mathrm{BF}$ such that for all $Y \in \mathrm{BF}$,
$$Y \in \mathrm{Pow}(V_{\omega+\mathbf{k}(n)}^{\langle \mathrm{BF}, \epsilon \rangle}) \textrm{ if and only if } \langle \mathbf{BF}, \epsilon \rangle \models Y \in X.$$
Now, $\langle \mathbf{BF}, \epsilon \rangle \models X= V_{\omega+\mathbf{k}(n+1)}$.
\Square
\end{proof}

Our next task is to confirm that, in $\mathrm{NFU}$, an ordinal $\alpha$ is finite if and only if $\mathbf{k}(\alpha)$ is finite in the structure $\langle \mathrm{BF}, \epsilon \rangle$. 

\begin{Lemma1} \label{Th:FiniteOrdinalsInNFUAreFiniteOrdinalsInBF}
The theory $\mathrm{NFU}$ proves that for all ordinals $\alpha$,
$$(\alpha \textrm{ is finite}) \textrm{ if and only if }  \langle \mathrm{BF}, \epsilon \rangle \models (\mathbf{k}(\alpha) \textrm{ is finite}).$$
\end{Lemma1}

\begin{proof}
We work in $\mathrm{NFU}$. We prove the contra-positive of both directions. Let $\alpha$ be an ordinal. Let $X= \{ \beta \in \mathrm{ON} \mid \beta < \alpha \}$.\\
Suppose that $\alpha$ is infinite. Therefore $T^2(\alpha)$ is infinite.  Let $f: X \longrightarrow X$ be an injection that is not a surjection. Define $G \subseteq \mathrm{BF}$ such that for all $x \in \mathrm{BF}$,
$$\begin{array}{lcl}
x \in G & \textrm{ if and only if } & \langle \mathrm{BF}, \epsilon \rangle \models (x= \langle \mathbf{k}(\beta), \mathbf{k}(\gamma) \rangle) \textrm{ and } f(\beta)= \gamma.  
\end{array}$$
Theorem \ref{Th:SetOfBFEXTsLessThanT2VAreExtensions} implies that there exists $Y \in \mathrm{BF}$ such that for all $x \in \mathrm{BF}$,
$$x \in G \textrm{ if and only if } \langle \mathrm{BF}, \epsilon \rangle \models x \in Y.$$
Now,
$$\langle \mathrm{BF}, \epsilon \rangle \models (Y \textrm{ is an injection from } \mathbf{k}(\alpha) \textrm{ into } \mathbf{k}(\alpha) \textrm{ that is not surjective}).$$
Conversely, suppose that $\alpha$ is such that
$$\langle \mathrm{BF}, \epsilon \rangle \models (\mathbf{k}(\alpha) \textrm{ is infinite}).$$
Let $G \in \mathrm{BF}$ be such that
$$\langle \mathrm{BF}, \epsilon \rangle \models (G \textrm{ is an injection from } \mathbf{k}(\alpha) \textrm{ into } \mathbf{k}(\alpha) \textrm{ that is not surjective}).$$
Define $f: X \longrightarrow X$ such that for all $\beta, \gamma < \alpha$,
$$f(\beta)= \gamma \textrm{ if and only if } \langle \mathrm{BF}, \epsilon \rangle \models (G(\mathbf{k}(\beta))= \mathbf{k}(\gamma)).$$
Now, $f$ witnesses the fact that $T^2(\alpha)$ is infinite. Therefore $\alpha$ is infinite.  
\Square
\end{proof}  

\noindent Theorem \ref{Th:AxCountLEQImpliesVOmegaPlusNExists} and Lemma \ref{Th:FiniteOrdinalsInNFUAreFiniteOrdinalsInBF} imply Hinnion's result about the strength of $\mathrm{AxCount}_\leq$ in the context of $\mathrm{NFU}$. 

\begin{Theorems1} \label{Th:AxCountLEQProveVOmegaExists}
(Hinnion \cite{hin75}) $$\mathrm{NFU}+ \mathrm{AxCount}_\leq \vdash (\langle \mathrm{BF}, \epsilon \rangle \models V_{\omega+\omega} \textrm{ exists}).$$
\end{Theorems1}

\begin{proof}
We work in $\mathrm{NFU}+\mathrm{AxCount}_\leq$. Theorem \ref{Th:AxCountLEQImpliesVOmegaPlusNExists} combined with Lemma \ref{Th:FiniteOrdinalsInNFUAreFiniteOrdinalsInBF} imply that 
$$\langle \mathrm{BF}, \epsilon \rangle \models (\forall n \in \mathbb{N})(V_{\omega+n} \textrm{ exists}).$$
Therefore, by collection in $\langle \mathrm{BF}, \epsilon \rangle$, 
$$\langle \mathrm{BF}, \epsilon \rangle \models V_{\omega+\omega} \textrm{ exists}.$$
\Square
\end{proof}

\begin{Coroll1} \label{Th:NFUPlusAxCountLEQProvesConZ}
(Hinnion \cite{hin75}) $$\mathrm{NFU}+\mathrm{AxCount}_\leq \vdash \mathrm{Con}(\mathrm{Z}).$$
\Square
\end{Coroll1}

\noindent It follows from Theorem \ref{Th:AutomorphismMovingNoNaturalDown} that the structure $\langle \mathrm{BF}, \epsilon\rangle$ in the theory $\mathrm{NFU}+\mathrm{AxCount}_\leq$ can also see a model of $\mathrm{NFU}+ \mathrm{AxCount}_\geq$.

\begin{Theorems1}
$$\mathrm{NFU}+\mathrm{AxCount}_\leq \vdash \mathrm{Con}(\mathrm{NFU}+\mathrm{AxCount}_\geq).$$
\end{Theorems1}

\begin{proof}
We work inside $\langle \mathrm{BF}, \epsilon \rangle$. Theorem \ref{Th:AxCountLEQProveVOmegaExists} shows that $V_{\omega+\omega}$ exists. Using the same arguments used in proof of Theorem \ref{Th:ModelOfNFUAxCountGEQBethOmegaDoesNotExist} we can see that this implies that there exists an $\mathcal{L}_{\mathcal{S}, P}$-structure $\mathcal{M}$ such that 
$$\mathcal{M} \models \mathrm{NFU}+ \mathrm{AxCount}_\geq.$$
\Square
\end{proof}

\noindent These results still leave the following questions unanswered:

\begin{Quest1}
What is the exact strength of the theory $\mathrm{NFU}+\mathrm{AxCount}_\geq$ relative to a subsystem of $\mathrm{ZFC}$?
\end{Quest1}

\begin{Quest1}
Does $\mathrm{NFU}+\mathrm{AxCount}_\geq$ prove the consistency of $\mathrm{NFU}$?
\end{Quest1}

\section[The strength of $\mathrm{NFU}+\mathrm{AxCount}_\leq$]{The strength of $\mathrm{NFU}+\mathrm{AxCount}_\leq$} \label{Sec:NFUPlusAxCountLEQ}

The aim of this section is to use the results of section \ref{Sec:AutomorphismsFromNonStandardModels} to shed light on the strength of the theory $\mathrm{NFU}+\mathrm{AxCount}_\leq$. We begin by showing that there is a model $\mathrm{NFU}+\mathrm{AxCount}_\leq$ in which the order-type of $\mathrm{CNI}$ equipped with the natural ordering is recursive. We will then show that $\mathrm{NFU}+\mathrm{AxCount}$ proves the consistency of the theory $\mathrm{NFU}+\mathrm{AxCount}_\leq$.\\
\\
\indent We begin by introducing an extension of Zermelo set theory that includes function symbols whose intended interpretations are the rank function and the function $\alpha \mapsto V_\alpha$. This theory will have the property that the standard part of a non-standard $\omega$-model of this theory will be a model $\mathrm{KP}^{\mathcal{P}}$. Let $\hat{\rho}$ and $\bar{V}$ be new unary function symbols.

\begin{Definitions1}
We define Zermelo with ranks set theory ($\mathrm{ZR}$) to be the $\mathcal{L}_{\hat{\rho}, \bar{V}}$-theory with axioms
\begin{itemize}
\item[(i)] all of the axioms of $\mathrm{Z}$,  
\item[(ii)] $\forall x(\hat{\rho}(x) \textrm{ is the least ordinal s.t. } (\forall y \in x)(\hat{\rho}(y) < \hat{\rho}(x)))$,
\item[(iii)] $\forall x((x \textrm{ is an ordinal}) \Rightarrow \forall y(y \in \bar{V}(x) \iff \hat{\rho}(y) < x))$,
\item[(iv)] separation for all $\mathcal{L}_{\hat{\rho}, \bar{V}}$-formulae.
\end{itemize}  
\end{Definitions1}

\noindent By considering the set $V_{\omega+\omega}$ and interpreting $\hat{\rho}$ using $\rho$ and $\bar{V}$ using the function $\alpha \mapsto V_\alpha$ we see that $\mathrm{KP}^\mathcal{P}$ proves the consistency of the theory $\mathrm{ZR}$.\\
\\
\indent We are now in a position to define the standard part of a non-standard model of $\mathrm{ZR}$.

\begin{Definitions1} \label{Df:StandardOrdinal}
Let $\mathcal{M}$ be an $\mathcal{L}_{\hat{\rho}, \bar{V}}$-structure such that $\mathcal{M} \models \mathrm{ZR}$. We define
$$\varpi_o(\mathcal{M})= \{\alpha \mid (\exists x \in M^\mathcal{M})((\mathcal{M} \models x \textrm{ is an ordinal}) \land (\langle x, \in^\mathcal{M} \rangle \cong \langle \alpha, \in \rangle)) \}.$$
\end{Definitions1}

\begin{Definitions1} \label{Df:StandardPartZR}
Let $\mathcal{M}$ be an $\mathcal{L}_{\hat{\rho}, \bar{V}}$-structure such that $\mathcal{M} \models \mathrm{ZR}$. We define
$$\mathbf{std}^{\mathrm{ZR}}(\mathcal{M})= \{ x \in M^\mathcal{M} \mid (\exists \alpha \in \varpi_o(\mathcal{M}))(\langle \hat{\rho}^\mathcal{M}(x), \in^\mathcal{M} \rangle \cong \langle \alpha, \in \rangle)\}.$$ 
\end{Definitions1}

We are now able to prove that the standard part of a non-standard $\omega$-model of $\mathrm{ZR}$ is a model of $\mathrm{KP}^{\mathcal{P}}$. The idea of building models of $\mathrm{KP}^{\mathcal{P}}$ from the standard part of a non-standard model first appears in \cite{fri73}.

\begin{Lemma1} \label{Th:StandardPartOfModelOfZR}
If $\mathcal{M} \models \mathrm{ZR}$ is a non-standard $\omega$-model then
$$\langle \mathbf{std}^{\mathrm{ZR}}(\mathcal{M}), \in^\mathcal{M} \rangle \models \mathrm{KP}^\mathcal{P}.$$
\end{Lemma1}

\begin{proof}
Let $\mathcal{M}= \langle M^\mathcal{M}, \in^\mathcal{M}, \hat{\rho}^\mathcal{M}, \bar{V}^\mathcal{M} \rangle$ be a non-standard $\omega$-model such that $\mathcal{M} \models \mathrm{ZR}$. Let $\mathcal{N}= \langle \mathbf{std}^{\mathrm{ZR}}(\mathcal{M}), \in^\mathcal{M} \rangle$. From the point of view of $\mathcal{M}$, $\mathbf{std}^{\mathrm{ZR}}(\mathcal{M})$ is a transitive subclass of the universe such that for all $\alpha \in \mathbf{std}^{\mathrm{ZR}}(\mathcal{M})$ with $\mathcal{N} \models \alpha \textrm{ is an ordinal}$, 
$$\bar{V}^\mathcal{M}(\alpha) \in \mathbf{std}^{\mathrm{ZR}}(\mathcal{M}).$$ 
Therefore $\mathcal{N} \models \mathrm{Mac}$ and if $\phi(\vec{x})$ is a $\Delta_0^\mathcal{P}$-formula then for all $\vec{a} \in \mathbf{std}^{\mathrm{ZR}}(\mathcal{M})$,
$$\mathcal{N} \models \phi(\vec{a}) \textrm{ if and only if } \mathcal{M} \models \phi(\vec{a}).$$
And so, $\mathcal{N}$ also satisfies $\Delta_0^\mathcal{P}$-Separation. We need to verify that $\Delta_0^\mathcal{P}$-Collection holds in $\mathcal{N}$. Let $\phi(x, y, \vec{z})$ be a $\Delta_0^\mathcal{P}$-formula. Let $\vec{a}, b \in \mathbf{std}^{\mathrm{ZR}}(\mathcal{M})$ be such that
\begin{equation} \label{Df:CollectionInStandardModelOfZR}
\mathcal{N} \models (\forall x \in b) \exists y \phi(x, y, \vec{a}).
\end{equation}
Let $\alpha \in M^\mathcal{M}$ be non-standard such that $\mathcal{M} \models \alpha \textrm{ is an ordinal}$. Working inside $\mathcal{M}$, let
$$A= \{ \beta \in \bar{V}(\alpha) \mid (\exists x \in b)(\exists y(\phi(x, y, \vec{a}) \land \hat{\rho}(y)= \beta) \land \forall y(\phi(x, y, \vec{a}) \Rightarrow \hat{\rho}(y) \geq \beta)) \}$$
$$\textrm{and } B= \{ \beta \in \bar{V}(\alpha) \mid (\beta \textrm{ is an ordinal}) \land (\beta \notin A) \}.$$
The fact that $\mathcal{M} \models \mathrm{ZR}$ ensures that $A$ and $B$ are sets. The fact that $\phi$ is absolute between $\mathcal{M}$ and $\mathcal{N}$ and (\ref{Df:CollectionInStandardModelOfZR}) ensures that $A^* \subseteq \mathbf{std}^{\mathrm{ZR}}(\mathcal{M})$. Let $\gamma$ be the $\in^\mathcal{M}$-least element of $B$. Therefore $\gamma \in \mathbf{std}^{\mathrm{ZR}}(\mathcal{M})$. Therefore, working inside $\mathcal{N}$, let 
$$C= \{ y \in \bar{V}(\gamma) \mid (\exists x \in b) \phi(x, y, \vec{a}) \}.$$
$\Delta_0^\mathcal{P}$-separation in $\mathcal{N}$ ensures that $C$ is a set in $\mathcal{N}$. Moreover,
$$\mathcal{N} \models (\forall x \in b)(\exists y \in C) \phi(x, y, \vec{a}).$$
Therefore $\Delta_0^{\mathcal{P}}$-collection holds in $\mathcal{N}$.\\
The fact that the relation $\in^\mathcal{N}$ is well-founded on $\mathbf{std}^{\mathrm{ZR}}(\mathcal{M})$ guarantees that $\mathcal{N}$ satisfies full class foundation.        
\Square
\end{proof}

Using the Barwise Compactness Theorem we can prove (in $\mathrm{ZFC}^-$) that if $V_{\omega_1^{\mathrm{ck}}}$ exists then there exists a non-standard $\omega$-model of $\mathrm{ZR}$.

\begin{Lemma1} \label{Th:NonStandardOmegaModelInZFCMinus}
The theory $\mathrm{ZFC}^-+ V_{\omega_1^{\mathrm{ck}}} \textrm{ exists}$ proves that there exists an $\mathcal{L}_{\hat{\rho}, \bar{V}}$-structure $\mathcal{N}$ such that $\mathcal{N} \models \mathrm{ZR}$ is a non-standard $\omega$-model. 
\end{Lemma1}

\begin{proof}
We work in the theory $\mathrm{ZFC}^-+ V_{\omega_1^{\mathrm{ck}}} \textrm{ exists}$. Let $\rho$ be the rank function. Define $V_*: \omega_1^{\mathrm{ck}} \longrightarrow V_{\omega_1^{\mathrm{ck}}}$ by $\alpha \mapsto V_\alpha$. The $\mathcal{L}_{\hat{\rho}, \bar{V}}$-structure $\mathcal{M}= \langle V_{\omega_1^{\mathrm{ck}}}, \in, \rho \cap (V_{\omega_1^{\mathrm{ck}}} \times V_{\omega_1^{\mathrm{ck}}}), V_* \rangle$ is such that $\mathcal{M} \models \mathrm{ZR}$.\\ 
Let $A= L_{\omega_1^{\mathrm{ck}}}$. Let $\mathcal{L}^\prime$ be the extension of $\mathcal{L}_{\hat{\rho}, \bar{V}}$ obtained by adding
\begin{itemize}
\item[(i)] constant symbols $\hat{\alpha}$ for every $\alpha \in \omega_1^{\mathrm{ck}}$,
\item[(ii)] constant symbol $\mathbf{c}$.
\end{itemize}
Define the $(\mathcal{L}^\prime_{\omega_1 \omega})_A$-theory $Q$ with axioms
\begin{itemize}
\item[(i)] all of the axioms of $\mathrm{ZR}$, 
\item[(ii)] $\hat{\zeta} \in \hat{\nu}$ for all $\zeta \in \nu \in \omega_1^{\mathrm{ck}}$, 
\item[(iii)] for all $\nu \in \omega_1^{\mathrm{ck}}$,
$$\forall x \left( x \in \hat{\nu} \Rightarrow \bigvee_{\zeta \in \nu} (x= \hat{\zeta})\right),$$
\item[(iv)] $\mathbf{c}$ is an ordinal,
\item[(v)] $\hat{\nu} \in \mathbf{c}$ for all $\nu \in \omega_1^{\mathrm{ck}}$,
\item[(vi)] every ordinal is recursive. 
\end{itemize}
The theory $Q$ is $\Sigma_1(\mathcal{L})$ over $A$. Let $Q^\prime \subseteq Q$ be such that $Q^\prime \in A$. Therefore the set
$$B= \{ \alpha \in \omega_1^{\mathrm{ck}} \mid \hat{\alpha} \textrm{ is mentioned in } Q^\prime \}$$
is bounded in $\omega_1^{\mathrm{ck}}$. Therefore we can expand $\mathcal{M}$ to a structure $\mathcal{M}^\prime$ that satisfies $Q^\prime$. Therefore by the Barwise Compactness Theorem (Theorem \ref{Th:BarwiseCompactness}) there exists a structure $\mathcal{N}^\prime$ that satisfies $Q$. Let $\mathcal{N}= \langle N^\mathcal{N}, \in^\mathcal{N}, \hat{\rho}^\mathcal{N}, \bar{V}^\mathcal{N} \rangle$ be the $\mathcal{L}_{\hat{\rho}, \bar{V}}$ reduct of $\mathcal{N}^\prime$. It follows from axiom scheme (i) of $Q$ that $\mathcal{N} \models \mathrm{ZR}$. Axiom schemes (ii) and (iii) of $Q$ ensure that $\mathcal{N}$ is an $\omega$-model. Let
$$C= \{ \beta \in N^\mathcal{N} \mid (\mathcal{N}^\prime \models \beta \textrm{ is an ordinal})\land (\forall \nu \in \omega_1^{\mathrm{ck}}) (\mathcal{N}^\prime \models \hat{\nu} \in \beta) \}.$$
By (iv) and scheme (v) of $Q$, $C$ is non-empty. Suppose that $\mathcal{N}$ is well-founded. Let $\alpha \in N^\mathcal{N}$ be the least element of $C$. But, $\langle \alpha, \in^\mathcal{N} \rangle \cong \langle \omega_1^{\mathrm{ck}}, \in \rangle$, which contradicts (vi) of $Q$. Therefore $\mathcal{N}$ is a non-standard $\omega$-model.   
\Square
\end{proof}

\noindent Combining Lemma \ref{Th:StandardPartOfModelOfZR} and Lemma \ref{Th:NonStandardOmegaModelInZFCMinus} immediately yields:

\begin{Coroll1}
$\mathrm{ZFC}^-+ V_{\omega_1^{\mathrm{ck}}} \textrm{ exists} \vdash \mathrm{Con}(\mathrm{KP}^{\mathcal{P}})$. \Square
\end{Coroll1}

The non-standard $\omega$-model built in Lemma \ref{Th:NonStandardOmegaModelInZFCMinus} also has the property that every ordinal is recursive. By assuming that we are working in a model of the generalised continuum hypothesis we can also ensure that every set of infinite sets with distinct cardinalities is isomorphic to a recursive ordinal when ordered by the natural ordering on cardinals.   

\begin{Lemma1} \label{Th:NonStandardOmegaModelFromZFCPlusGCH}
The theory $\mathrm{ZFC}+\mathrm{GCH}$ proves that there exists an $\mathcal{L}_{\hat{\rho}, \bar{V}}$-structure $\mathcal{N}$ such that 
$$\mathcal{N} \models \mathrm{ZR} +\textrm{every ordinal is recursive},$$
\begin{equation} \label{Df:CNIAreRecursive}
\textrm{and } \mathcal{N} \models \forall X \left( \begin{array}{c}
(\forall x, y \in X)(|x| \geq \omega \land |y| \geq \omega \land (|x| < |y| \lor |y| < |x| \lor x=y)) \Rightarrow\\
\exists f \exists \alpha\left( \begin{array}{c}
(\alpha \textrm{ is a recursive ordinal}) \land (f: X \longrightarrow \alpha)\\
\land (\forall x, y \in X)(|x| < |y| \Rightarrow f(x) < f(y))
\end{array}\right) 
\end{array}\right),
\end{equation}
and $\mathcal{N}$ is a non-standard $\omega$-model.
\end{Lemma1}

\begin{proof}
We work in the theory $\mathrm{ZFC}+\mathrm{GCH}$. Let $\rho$ be the rank function. Define $V_*: \omega_1^{\mathrm{ck}} \longrightarrow V_{\omega_1^{\mathrm{ck}}}$ by $\alpha \mapsto V_\alpha$. The $\mathcal{L}_{\hat{\rho}, \bar{V}}$-structure $\mathcal{M}= \langle V_{\omega_1^{\mathrm{ck}}}, \in, \rho \cap (V_{\omega_1^{\mathrm{ck}}} \times V_{\omega_1^{\mathrm{ck}}}), V_* \rangle$ is such that $\mathcal{M} \models \mathrm{ZR}$ and $\mathcal{M} \models \textrm{every ordinal is recursive}$. The Generalised Continuum Hypothesis ensures that $\mathcal{M}$ satisfies (\ref{Df:CNIAreRecursive}). Using the same arguments used in the proof of Lemma \ref{Th:NonStandardOmegaModelInZFCMinus} we can find an $\mathcal{L}_{\hat{\rho}, \bar{V}}$-structure $\mathcal{N}$ such that 
$$\mathcal{N} \models \mathrm{ZR} +\textrm{every ordinal is recursive},$$
$$
\textrm{and } \mathcal{N} \models \forall X \left( \begin{array}{c}
(\forall x, y \in X)(|x| \geq \omega \land |y| \geq \omega \land (|x| < |y| \lor |y| < |x| \lor x=y)) \Rightarrow\\
\exists f \exists \alpha\left( \begin{array}{c}
(\alpha \textrm{ is a recursive ordinal}) \land (f: X \longrightarrow \alpha)\\
\land (\forall x, y \in X)(|x| < |y| \Rightarrow f(x) < f(y))
\end{array}\right) 
\end{array}\right),
$$
and $\mathcal{N}$ is a non-standard $\omega$-model.
\Square
\end{proof}

\noindent By applying Theorem \ref{Th:AutomorphismsFromNonStandardOmegaModels} to the model $\mathcal{N}$ built in Lemma \ref{Th:NonStandardOmegaModelFromZFCPlusGCH} we obtain a model of $\mathrm{NFU}+\mathrm{AxCount}_\leq$ in which $\langle \mathrm{CNI}, < \rangle$ is isomorphic to a recursive ordinal. 

\begin{Theorems1} \label{Th:WeakModelOfNFUPlusAxCountLEQ}
There exists an $\mathcal{L}_{\mathcal{S}, P}$-structure $\mathcal{M}$ such that $\mathcal{M} \models \mathrm{NFU}+\mathrm{AxCount}_\leq$ and
$$\mathcal{M} \models \langle \mathrm{CNI}, < \rangle \textrm{ is isomorphic to a recursive ordinal}.$$
\end{Theorems1}

\begin{proof}
We work in the theory $\mathrm{ZFC}+\mathrm{GCH}$. By Lemma \ref{Th:NonStandardOmegaModelFromZFCPlusGCH} there exists a non-standard $\omega$-model $\mathcal{Q} \models \mathrm{ZR} + \textrm{every ordinal is recursive}$ such that
\begin{equation} \label{Df:CNIRecursiveII}
\mathcal{Q} \models \forall X \left( \begin{array}{c}
(\forall x, y \in X)(|x| \geq \omega \land |y| \geq \omega \land (|x| < |y| \lor |y| < |x| \lor x=y)) \Rightarrow\\
\exists f \exists \alpha\left( \begin{array}{c}
(\alpha \textrm{ is a recursive ordinal}) \land (f: X \longrightarrow \alpha)\\
\land (\forall x, y \in X)(|x| < |y| \Rightarrow f(x) < f(y))
\end{array}\right) 
\end{array}\right),
\end{equation}
Let $\mathcal{Q}^\prime$ be the $\mathcal{L}$ reduct of $\mathcal{Q}$. Let $\mathcal{Q}^{\prime\prime}$ be the expansion of $\mathcal{Q}^\prime$ to an $\mathcal{L}_{\bar{\beta}}$-structure satisfying the assumptions of Theorem \ref{Th:AutomorphismsFromNonStandardOmegaModels}. Therefore by Theorem \ref{Th:AutomorphismsFromNonStandardOmegaModels} there exists an $\mathcal{L}_{\bar{\beta}}$-structure $\mathcal{N} \equiv \mathcal{Q}^{\prime\prime}$ admitting an automorphism $j: \mathcal{N} \longrightarrow \mathcal{N}$ such that
\begin{itemize}
\item[(i)] $\mathcal{N} \models j(n) \geq n$ for all $n \in (\omega^\mathcal{N})^*$, 
\item[(ii)] there exists $\alpha \in (\bar{\beta}^\mathcal{N})^*$ such that $\mathcal{N} \models j(\alpha) < \alpha$. 
\end{itemize}
Since $\mathcal{N} \equiv \mathcal{Q}^{\prime\prime}$, the rank function $\rho$ is well-defined and 
$$\mathcal{N} \models \forall \alpha ((\alpha \textrm{ is an ordinal}) \Rightarrow (V_\alpha \textrm{ exists})).$$
Let $\alpha \in (\bar{\beta}^\mathcal{N})^*$ be such that $\mathcal{N} \models j(\alpha) < \alpha$. Recalling (\ref{Df:TZTUFromMacWithAut1}) and (\ref{Df:TZTUFromMacWithAut2}) we can use $j^{-1}$ and $V_{\alpha}^\mathcal{N}$ to define an $\mathcal{L}_{\mathrm{T}\mathbb{Z}\mathrm{TU}}$-structure $\mathcal{M}_{j^{-1}}^{V_\alpha^\mathcal{N}} \models \mathrm{T}\mathbb{Z}\mathrm{TU}$ such that $j^{-1}$ is a type-shifting automorphism of $\mathcal{M}_{j^{-1}}^{V_\alpha^\mathcal{N}}$. Using (\ref{Df:NFUMembership}) and Theorem \ref{Th:NFUFromMacWithAut} we can define $\in_{j^{-1}}$, $S$ and $\bar{P}$ on $V_{\alpha}^\mathcal{N}$ so as 
$$\mathcal{M}= \langle V_{\alpha}^\mathcal{N}, \in_{j^{-1}}, S, \bar{P} \rangle \models \mathrm{NFU}.$$
By the properties of $j^{-1}$, for all $n \in (\omega^\mathcal{N})^*$, 
$$\mathcal{M}^{V_{\alpha}^{\mathcal{N}}}_{j^{-1}} \models j^{-1}(n) \leq |\iota``n|.$$
Therefore, by Theorem \ref{Th:NFUFromMacWithAut}, for all $n \in (\omega^\mathcal{N})^*$,
$$\mathcal{M} \models T(n) \geq n.$$
Therefore $\mathcal{M} \models \mathrm{AxCount}_\leq$. Now, let $\sigma_i$ by the $\mathcal{L}_{\mathrm{T}\mathbb{Z}\mathrm{TU}}$-sentence that says: `the set of infinite cardinals [at level $i$] equipped with the natural ordering is isomorphic to a recursive ordinal'. It follows from (\ref{Df:CNIRecursiveII}) that for all $i \in \mathbb{Z}$,
$$\mathcal{M}^{V_{\alpha}^{\mathcal{N}}}_{j^{-1}} \models \sigma_i.$$
Therefore 
$$\mathcal{M} \models \textrm{the order-type of } \langle \mathrm{CNI}, < \rangle \textrm{ is recursive}.$$ 
\Square
\end{proof}

\noindent Again, since $\bar{\bar{\beth}} \subseteq \mathrm{CNI}$, we get:

\begin{Coroll1}
There exists an $\mathcal{L}_{\mathcal{S}, P}$-structure $\mathcal{M}$ such that $\mathcal{M} \models \mathrm{NFU}+\mathrm{AxCount}_\leq$ and
$$\mathcal{M} \models \beth_{\omega_1^{\mathrm{ck}}}^{TT} \textrm{ does not exist}.$$
\Square
\end{Coroll1}

Our aim now turns to proving that the theory $\mathrm{NFU}+\mathrm{AxCount}$ is strictly stronger than $\mathrm{NFU}+\mathrm{AxCount}_\leq$. Our tactic for proving this result will be to apply Lemma \ref{Th:NonStandardOmegaModelInZFCMinus} and Theorem \ref{Th:AutomorphismsFromNonStandardOmegaModels} inside the structure $\langle \mathrm{BF}, \epsilon \rangle$ in the theory $\mathrm{NFU}+\mathrm{AxCount}$. In order to apply Lemma \ref{Th:NonStandardOmegaModelInZFCMinus} inside $\langle \mathrm{BF}, \epsilon \rangle$ we first need to show that $\mathrm{NFU}+\mathrm{AxCount}$ proves that 
$$\langle \mathrm{BF}, \epsilon \rangle \models V_{\omega_1^{\mathrm{ck}}} \textrm{ exists}.$$
We will prove this result using the same tactics employed in section \ref{Sec:NFUPlusAxCountGEQ}.

\begin{Lemma1} \label{Th:CountableOrdinalsAreStronglyCantorian}
The theory $\mathrm{NFU}+\mathrm{AxCount}$ proves that if $\alpha$ is a countable then\\
$T(\alpha)= \alpha$.
\end{Lemma1}

\begin{proof}
We work in the theory $\mathrm{NFU}+\mathrm{AxCount}$. Let $\alpha$ be a countable ordinal. Let $R \in \alpha$ be such that $\mathrm{Dom}(R) \subseteq \mathbb{N}$. Let 
$$S= \{ \langle \{x\}, \{y\} \rangle \mid \langle x, y \rangle \in R \}.$$
Define $f: \mathrm{Dom}(R) \longrightarrow \mathrm{Dom}(S)$ by
$$\begin{array}{lcl}
f(x)= \{y\} & \textrm{ if and only if } & T(y)=x
\end{array}.$$
Stratified Comprehension ensures that $f$ exists. Since $T$ is the identity on $\mathrm{Dom}(R)$, $f$ witnesses the fact that $R \cong S$. Therefore $T(\alpha)= \alpha$.
\Square
\end{proof}

\noindent This allows us to prove an extension of Lemma \ref{Th:BethNLessThanT2V} in the presence of $\mathrm{AxCount}$.

\begin{Lemma1} \label{Th:BethAlphaLEQT2V}
The theory $\mathrm{NFU}+\mathrm{AxCount}$ proves that for all countable ordinals $\alpha$, $\beth_\alpha^{TT} \leq T^2(|V|)$.
\end{Lemma1}

\begin{proof}
We work in the theory $\mathrm{NFU}+\mathrm{AxCount}$. Let $\alpha$ be the least countable ordinal such that $T^2(|V|) < \beth_\alpha^{TT} \leq |V|$. But this immediately leads to contradiction because, 
$$\beth_\alpha^{TT}= \beth_{T^2(\alpha)}^{TT}= T^2(\beth_\alpha^{TT}) \leq T^2(|V|).$$ 
\Square
\end{proof}

\noindent This extension facilitates the proof of an analogue of Theorem \ref{Th:AxCountLEQImpliesVOmegaPlusNExists}.

\begin{Theorems1} \label{Th:AxCountProvesCountableRanksExist}
The theory $\mathrm{NFU}+\mathrm{AxCount}$ proves that for all countable ordinals $\alpha$,
\begin{itemize}
\item[(i)] $\langle \mathrm{BF}, \epsilon \rangle \models V_{\mathbf{k}(\alpha)} \textrm{ exists}$,
\item[(ii)] $|(V_{\mathbf{k}(\alpha)}^{\langle \mathrm{BF}, \epsilon \rangle})^{\mathrm{ext}}|= \beth_\alpha^{TT}$. 
\end{itemize}
\end{Theorems1}

\begin{proof}
We work in the theory $\mathrm{NFU}+\mathrm{AxCount}$. We prove the theorem by induction. Corollary \ref{Th:BFEXTsModelZFCminus} implies that 
\begin{itemize}
\item[(i)] $\langle \mathrm{BF}, \epsilon \rangle \models V_\omega \textrm{ exists}$,
\item[(ii)] $|(V_{\omega}^{\langle \mathrm{BF}, \epsilon \rangle})^{\mathrm{ext}}|= \aleph_0= \beth_0^{TT}$.
\end{itemize}
Let $\alpha$ be a countable ordinal and suppose that the Theorem holds for all $\beta < \alpha$. If $\alpha$ is a limit ordinal then collection in $\langle \mathrm{BF}, \epsilon \rangle$ implies that
$$\langle \mathrm{BF}, \epsilon \rangle \models V_{\mathbf{k}(\alpha)} \textrm{ exists}.$$
Moreover, 
$$|(V_{\mathbf{k}(\alpha)}^{\langle \mathrm{BF}, \epsilon \rangle})^{\mathrm{ext}}|= \sup \{ |(V_{\mathbf{k}(\beta)}^{\langle \mathrm{BF}, \epsilon \rangle})^{\mathrm{ext}}| \mid \beta < \alpha \}= \beth_\alpha^{TT}.$$
Now, suppose that $\alpha= \gamma + 1$. By Lemma \ref{Th:SizeOfPow},
$$|\mathrm{Pow}(V_{\mathbf{k}(\gamma)}^{\langle \mathrm{BF}, \epsilon \rangle})| = 2^{|(V_{\mathbf{k}(\gamma)}^{\langle \mathrm{BF}, \epsilon \rangle})^{\mathrm{ext}}|}= 2^{\beth_\gamma^{TT}}= \beth_{\alpha}^{TT}.$$
Therefore, by Lemma \ref{Th:BethAlphaLEQT2V}, $|\mathrm{Pow}(V_{\mathbf{k}(\gamma)}^{\langle \mathrm{BF}, \epsilon \rangle})| \leq T^2(|V|)$. And so, by Theorem \ref{Th:SetOfBFEXTsLessThanT2VAreExtensions}, there exists an $X \in \mathrm{BF}$ such that for all $Y \in \mathrm{BF}$,
$$Y \in \mathrm{Pow}(V_{\mathbf{k}(\gamma)}^{\langle \mathrm{BF}, \epsilon \rangle}) \textrm{ if and only if } \langle \mathbf{BF}, \epsilon \rangle \models Y \in X.$$
Now, $\langle \mathbf{BF}, \epsilon \rangle \models X= V_{\mathbf{k}(\alpha)}$.
\Square
\end{proof}

We now turn to extending Lemma \ref{Th:FiniteOrdinalsInNFUAreFiniteOrdinalsInBF} to show that, in $\mathrm{NFU}$, an ordinal $\alpha$ is countable if and only if $\mathbf{k}(\alpha)$ is countable in $\langle \mathrm{BF}, \epsilon \rangle$.

\begin{Lemma1} \label{Th:CountableOrdinalsInNFUAreCountableInBF}
The theory $\mathrm{NFU}$ proves that for all ordinals $\alpha$,
$$(\alpha \textrm{ is countable}) \textrm{ if and only if } \langle \mathrm{BF}, \epsilon \rangle \models (\mathbf{k}(\alpha) \textrm{ is countable}).$$
\end{Lemma1}

\begin{proof}
We work in $\mathrm{NFU}$. Let $\alpha$ be an ordinal. Let $X= \{ \beta \in \mathrm{ON} \mid \beta < \alpha \}$.\\
Suppose that $\alpha$ is countable. Therefore $T^2(\alpha)$ is countable.  Let $f: X \longrightarrow \mathbb{N}$ be an injection. Define $G \subseteq \mathrm{BF}$ such that for all $x \in \mathrm{BF}$,
$$\begin{array}{lcl}
x \in G & \textrm{ if and only if } & \langle \mathrm{BF}, \epsilon \rangle \models (x= \langle \mathbf{k}(\beta), \mathbf{k}(n) \rangle) \textrm{ and } f(\beta)= n.  
\end{array}$$
Theorem \ref{Th:SetOfBFEXTsLessThanT2VAreExtensions} implies that there exists $Y \in \mathrm{BF}$ such that for all $x \in \mathrm{BF}$,
$$x \in G \textrm{ if and only if } \langle \mathrm{BF}, \epsilon \rangle \models x \in Y.$$
Now,
$$\langle \mathrm{BF}, \epsilon \rangle \models (Y \textrm{ is an injection from } \mathbf{k}(\alpha) \textrm{ into } \omega).$$
Conversely, suppose that $\alpha$ is such that
$$\langle \mathrm{BF}, \epsilon \rangle \models (\mathbf{k}(\alpha) \textrm{ is countable}).$$
Let $G \in \mathrm{BF}$ be such that
$$\langle \mathrm{BF}, \epsilon \rangle \models (G \textrm{ is an injection from } \mathbf{k}(\alpha) \textrm{ into } \omega).$$
Define $f: X \longrightarrow \mathbb{N}$ such that for all $\beta < \alpha$ and for all $n \in \mathbb{N}$,
$$f(\beta)= n \textrm{ if and only if } \langle \mathrm{BF}, \epsilon \rangle \models (G(\mathbf{k}(\beta))= \mathbf{k}(n)).$$
Now, Lemma \ref{Th:FiniteOrdinalsInNFUAreFiniteOrdinalsInBF} implies that $f$ is a function witnessing the fact that $T^2(\alpha)$ is countable. Therefore $\alpha$ is countable.
\Square
\end{proof}

\noindent Theorem \ref{Th:AxCountProvesCountableRanksExist} combined with Lemma \ref{Th:CountableOrdinalsInNFUAreCountableInBF} show that, in the theory $\mathrm{NFU}+\mathrm{AxCount}$, the structure $\langle \mathrm{BF}, \epsilon \rangle$ believes that $V_{\omega_1^{\mathrm{ck}}}$ exists.

\begin{Coroll1} \label{Th:NFUPlusAXCountProvesVOmega1CKExists}
$$\mathrm{NFU}+\mathrm{AxCount} \vdash (\langle \mathrm{BF}, \epsilon \rangle \models V_{\omega_1^{\mathrm{ck}}} \textrm{ exists}).$$
\Square 
\end{Coroll1}

\noindent Applying Lemma \ref{Th:StandardPartOfModelOfZR} inside the structure $\langle \mathrm{BF}, \epsilon \rangle$ immediately shows that adding $\mathrm{AxCount}$ to $\mathrm{NFU}$ allows us to strengthen Corollary \ref{Th:NFUPlusAxCountLEQProvesConZ}.

\begin{Coroll1}
$$\mathrm{NFU}+\mathrm{AxCount} \vdash \mathrm{Con}(\mathrm{KP}^\mathcal{P}).$$
\Square
\end{Coroll1}

\noindent Corollary \ref{Th:NFUPlusAXCountProvesVOmega1CKExists} also allows us to show that the theory $\mathrm{NFU}+\mathrm{AxCount}$ is strictly stronger than $\mathrm{NFU}+\mathrm{AxCount}_\leq$.

\begin{Theorems1}
$$\mathrm{NFU}+\mathrm{AxCount} \vdash \mathrm{Con}(\mathrm{NFU}+\mathrm{AxCount}_\leq).$$
\end{Theorems1}

\begin{proof}
We work in the theory $\mathrm{NFU}+\mathrm{AxCount}$. By Corollary \ref{Th:NFUPlusAXCountProvesVOmega1CKExists}:
$$\langle \mathrm{BF}, \epsilon \rangle \models V_{\omega_1^{\mathrm{ck}}} \textrm{ exists}.$$
By applying Lemma \ref{Th:NonStandardOmegaModelInZFCMinus} and Theorem \ref{Th:AutomorphismsFromNonStandardOmegaModels} and using the same arguments used in the proof of Theorem \ref{Th:WeakModelOfNFUPlusAxCountLEQ} we can build a model of $\mathrm{NFU}+\mathrm{AxCount}_\leq$ inside the structure $\langle \mathrm{BF}, \epsilon \rangle$.
\Square
\end{proof}

This section still leaves the following question unanswered:

\begin{Quest1} \label{Df:StrengthOfAxCountLEQQuestion}
What is the exact strength of the theory $\mathrm{NFU}+\mathrm{AxCount}_\leq$ relative to a subsystem of $\mathrm{ZFC}$?
\end{Quest1}

\noindent \textbf{Acknowledgements.} This research was completed while I was a Ph.D. student in the Department of Pure Mathematics and Mathematical Statistics at the University of Cambridge. I would like to thank my supervisor Thomas Forster for all of his help and support. I would also like to thank Ali Enayat, Randall Holmes and Andrey Bovykin for many helpful discussions. My Ph.D. studies were supported the Cambridge Commonwealth Trusts. 

\bibliographystyle{plain}
\bibliography{automorphismsandnonstandardmodels24}

\end{document}